\documentclass[10pt]{amsart}
\usepackage{amsfonts}
\usepackage{amsmath}

\newtheorem{theorem}{Theorem}
\newtheorem{definition}{Definition}
\newtheorem{proposition}{Proposition}

\begin{document}
\title{Symmetries of second order differential equations on Lie algebroids}
\author{Liviu Popescu}
\maketitle

\begin{abstract}
In this paper we investigate the relations between semispray, nonlinear
connection, dynamical covariant derivative and Jacobi endomorphism on Lie
algebroids. Using these geometric structures, we study the symmetries of
second order differential equations in the general framework of Lie
algebroids.
\end{abstract}



\section{\textbf{Introduction}}

The geometry of second order differential equations (SODE) on the
tangent bundle $TM$ of a differentiable manifold $M$ is closely
related to the geometry of nonlinear connections \cite{Cr1, Gr1}.
The system of SODE can by represented using the notion of
semispray, which together with the nonlinear connection induce two
important concepts: the dynamical covariant derivative and Jacobi
endomorphism \cite{Bu2, Bu3, Cr2, Gr2, Ma, Sa}. The notion of
dynamical covariant derivative was introduced for the first time
in the case of tangent bundle by J. Cari\~nena and E. Mart\'\i nez
\cite{Ca0} as a derivation of degree 0 along the tangent bundle
projection. The notion of symmetry in fields theory using various
geometric framework is intensely studied (see for instance
\cite{Ab, Bua, Bu3, Ga, Le0, Le1,Pr1, Pr2, Ro}). The notion of Lie
algebroid is a natural generalization of the tangent bundle and
Lie algebra. In the last decades the Lie algebroids \cite{Mk1,
Mk2} are the objects of intensive studies with applications to
mechanical systems or optimal control \cite {Ar2, Co, Fe, Le, Li,
Ma2, Ma4, Me, Po1, Po2, Po2', Po3, Po4, We} and are the natural
framework in which one can develop the theory of differential
equations, where the notion of symmetry plays a very important
role.

In this paper we study some properties of semispray and generalize
the notion of symmetry for second order differential equations on
Lie algebroids and characterize its properties using the dynamical
covariant derivative and Jacobi endomorphism. The paper is
organized as follows. In section two the preliminary geometric
structures on Lie algebroids are introduced and some relations
between them are given. We present the Jacobi endomorphism on Lie
algebroids and find the relation with the curvature tensor of
Ehresmann nonlinear connection. In section three we study the
dynamical covariant derivative on Lie algebroids. Using a
semispray and an arbitrary nonlinear connection, we introduce the
dynamical covariant derivative on Lie algebroids as a tensor
derivation and prove that the compatibility condition with the
tangent structure fixes the canonical nonlinear connection. In the
case of the canonical nonlinear connection induced by a semispray,
more properties of dynamical covariant derivative are added. In
the case of homogeneous second order differential equations
(spray) the relation between the dynamical covariant derivative
and Berwald linear connection is given. In the last section we
study the dynamical symmetries, Lie symmetries, Newtonoid sections
and Cartan symmetries on Lie algebroids and find the relations
between them. These structures are studied for the first time on
the tangent bundle by G. Prince in \cite{Pr1, Pr2}. Also, we prove
that an exact Cartan symmetry induces a conservation law and
conversely, which extends the work developed in \cite{Mar}.
Moreover, we find the invariant equations of dynamical symmetries,
Lie symmetries and Newtonoid sections in terms of dynamical
covariant derivative and Jacobi endomorphism, which generalize
some results from \cite{Bu3, Pr1, Pr2}. We have to mention that
the Noether type theorems for Lagrangian systems on Lie algebroids
can be found in \cite{Ca, Ma2} and Jacobi sections for second
order differential equations on Lie algebroids are studied in
\cite{Ca1}. Finally, using an example from optimal control theory
(driftless control affine systems), we prove that the framework of
Lie algebroids is more useful than the tangent bundle in order to
find the symmetries of the dynamics induced by a Lagrangian
function. Also, using the $k$-symplectic formalism on Lie
algebroids developed in \cite{Le2} one can study the symmetries in
this new framework, which generalize the results from \cite{Bua}.

\section{\textbf{Lie algebroids}}

Let $M$ be a real, $C^\infty $-differentiable, $n$-dimensional
manifold and $ (TM,\pi _M,M)$ its tangent bundle. A Lie algebroid
over a manifold $M$ is a triple $(E,[\cdot ,\cdot ]_E,\sigma )$,
where ($E,\pi ,M)$ is a vector bundle of rank $m$ over $M,$ which
satisfies the following conditions: \\a) $ C^\infty (M)$-module of
sections $\Gamma (E)$ is equipped with a Lie algebra structure
$[\cdot ,\cdot ]_E$. \\b) $\sigma :E\rightarrow TM$ is a bundle
map (called the anchor) which induces a Lie algebra homomorphism
(also denoted $\sigma $) from the Lie algebra of sections $(\Gamma
(E),[\cdot ,\cdot ]_E)$ to the Lie algebra of vector fields
$(\mathcal{\chi }(M),[\cdot ,\cdot ])$ satisfying the Leibnitz
rule
\begin{equation}
\lbrack s_1,fs_2]_E=f[s_1,s_2]_E+(\sigma (s_1)f)s_2,\ \forall s_1,s_2\in
\Gamma (E),\ f\in C^\infty (M).
\end{equation}

From the above definition it results: \\$1^{\circ }$ $[\cdot ,\cdot ]_E$ is
a $\Bbb{R}$-bilinear operation, \\$2^{\circ }$ $[\cdot ,\cdot ]_E$ is
skew-symmetric, i.e. $[s_1,s_2]_E=-[s_2,s_1]_E,\quad \forall s_1,s_2\in
\Gamma (E),$\\$3^{\circ }$ $[\cdot ,\cdot ]_E$ verifies the Jacobi identity
\begin{equation*}
\lbrack s_1,[s_2,s_3]_E]_E+[s_2,[s_3,s_1]_E]_E+[s_3,[s_1,s_2]_E]_E=0,\
\end{equation*}
and $\sigma $ being a Lie algebra homomorphism, means that $\sigma
[s_1,s_2]_E=[\sigma (s_1),\sigma (s_2)].$

The existence of a Lie bracket on the space of sections of a Lie
algebroid leads to a calculus on its sections analogous to the
usual Cartan calculus on differential forms. If $f$ is a function
on $M$, then $df(x)\in E_x^{*}$ is given by $ \left\langle
df(x),a\right\rangle =\sigma (a)f$, for $\forall a\in E_x$. For
$\omega $ $\in $ $\bigwedge^k(E^{*})$ the \textit{exterior
derivative} $ d^E\omega \in \bigwedge^{k+1}(E^{*})$ is given by
the formula
\begin{eqnarray*}
d^E\omega (s_1,...,s_{k+1}) &=&\overset{k+1}{\sum_{i=1}}(-1)^{i+1}\sigma
(s_i)\omega (s_1,...,\overset{\symbol{94}}{s}_i,...,s_{k+1})+ \\
&&\ \ \ \ \ \ \ +\sum_{1\leq i<j\leq k+1}(-1)^{i+j}\omega
([s_{i,}s_j]_E,s_1,...,\overset{\symbol{94}}{s_i},...,\overset{\symbol{94}}{
s_j},...s_{k+1}),
\end{eqnarray*}
where $s_i\in \Gamma (E)$, $i=\overline{1,k+1}$, and the hat over an
argument means the absence of the argument. It results that
\begin{equation*}
(d^E)^2=0,\quad d^E(\omega _1\wedge \omega _2)=d^E\omega _1\wedge
\omega _2+(-1)^{\deg \omega _1}\omega _1\wedge d^E\omega _2.
\end{equation*}
The cohomology associated with $d^E$ is called the \textit{Lie
algebroid cohomology} of $E$. Also, for $\xi $ $\in \Gamma (E)$
one can define the \textit{Lie derivative} with respect to $\xi$,
given by $\mathcal{L}_\xi =i_\xi \circ d^E+d^E\circ i_\xi $, where
$i_\xi $ is the contraction with $\xi $. We recall that if $L$ and
$K$ are $(1,1)$-type tensor field, Fr\"olicher-Nijenhuis bracket
$[L,K]$ is the vector valued 2-form \cite{Fr}
\begin{eqnarray*}
\lbrack L,K]_E(X,Y) &=&[LX,KY]_E+[KX,LY]_E+(LK+KL)[X,Y]_E- \\
&&\ \ \ \ \ \ -L[X,KY]_E-K[X,LY]_E-L[KX,Y]_E-K[LX,Y]_E,
\end{eqnarray*}
and the Nijenhuis tensor of $L$ is given by
\begin{equation*}
\mathbf{N}_L(X,Y)=\frac 12[L,L]_E=[LX,LY]_E+L^2[X,Y]_E-L[X,LY]_E-L[LX,Y]_E.
\end{equation*}
For a vector field in $\mathcal{X}(E)$ and a $(1,1)$-type tensor
field $L$ on $E$ the Fr\"olicher-Nijenhuis bracket
$[X,L]_E=\mathcal{L}_XL$ is the $ (1,1)$-type tensor field on $E$
given by
\begin{equation*}
\mathcal{L}_XL=\mathcal{L}_X\circ L-L\circ \mathcal{L}_X,
\end{equation*}
where $\mathcal{L}_X$ is the usual Lie derivative. If we take the
local coordinates $(x^i)$ on an open $U\subset $ $M$, a local
basis $\{s_\alpha \}$ of the sections of the bundle $\pi
^{-1}(U)\rightarrow U$ generates local coordinates $(x^i,y^\alpha
)$ on $E$. The local functions $\sigma _\alpha ^i(x)$, $L_{\alpha
\beta }^\gamma (x)$ on $M$ given by
\begin{equation}
\sigma (s_\alpha )=\sigma _\alpha ^i\frac \partial {\partial
x^i},\quad [s_\alpha ,s_\beta ]_E=L_{\alpha \beta }^\gamma
s_\gamma ,\quad i=\overline{1,n},\quad \alpha ,\beta ,\gamma
=\overline{1,m},
\end{equation}
are called the \textit{structure functions of the Lie algebroid},
and satisfy the \textit{structure equations} on Lie algebroids
\begin{equation*}
\sum_{(\alpha ,\beta ,\gamma )}\left( \sigma _\alpha ^i\frac{\partial
L_{\beta \gamma }^\delta }{\partial x^i}+L_{\alpha \eta }^\delta L_{\beta
\gamma }^\eta \right) =0,\quad \sigma _\alpha ^j\frac{\partial \sigma _\beta
^i}{\partial x^j}-\sigma _\beta ^j\frac{\partial \sigma _\alpha ^i}{\partial
x^j}=\sigma _\gamma ^iL_{\alpha \beta }^\gamma .
\end{equation*}
Locally, if $f\in C^\infty (M)$ then $d^Ef=\frac{\partial
f}{\partial x^i} \sigma _\alpha ^is^\alpha ,$ where $\{s^\alpha
\}$ is the dual basis of $ \{s_\alpha \}$ and if $\theta \in
\Gamma (E^{*}),$ $\theta =\theta _\alpha s^\alpha $ then
\begin{equation*}
d^E\theta =\left( \sigma _\alpha ^i\frac{\partial \theta _\beta
}{\partial x^i}-\frac 12\theta _\gamma L_{\alpha \beta }^\gamma
\right) s^\alpha \wedge s^\beta.
\end{equation*}
Particularly, we get $d^Ex^i=\sigma _\alpha ^is^\alpha $ and $d^Es^\alpha
=-\frac 12L_{\beta \gamma }^\alpha s^\beta \wedge s^\gamma .$

\subsection{\textbf{The prolongation of a Lie algebroid over the vector
bundle projection}}

Let $(E,\pi ,M)$ be a vector bundle. For the projection $\pi :E\rightarrow M$
we can construct the prolongation of $E$ (see \cite{Hi, Le, Ma2}). The
associated vector bundle is ($\mathcal{T}E,\pi _2,E$) where
\begin{equation*}
\mathcal{T}E=\underset{_{w\in E}}{\cup }\mathcal{T}_wE,\quad
\mathcal{T} _wE=\{(u_x,v_w)\in E_x\times T_wE\mid \sigma
(u_x)=T_w\pi (v_w),\quad \pi (w)=x\in M\},
\end{equation*}
and the projection $\pi _2(u_x,v_w)=\pi _E(v_w)=w$, where $\pi
_E:TE\rightarrow E$ is the tangent projection. We also have the
canonical projection $\pi _1:\mathcal{T}E\rightarrow E$ given by
$\pi _1(u,v)=u$. The projection onto the second factor $\sigma
^1:\mathcal{T}E\rightarrow TE$, $ \sigma ^1(u,v)=v$ will be the
anchor of a new Lie algebroid over the manifold $E$. An element of
$\mathcal{T}E$ is said to be vertical if it is in the kernel of
the projection $\pi _1$. We will denote $(V\mathcal{T}E,\pi
_{2\mid _{V\mathcal{T}E}},E)$ the vertical bundle of
$(\mathcal{T}E,\pi _2,E) $ and $\sigma ^1\left|
_{V\mathcal{T}E}\right. :V\mathcal{T} E\rightarrow VTE$ is an
isomorphism. If $f\in C^\infty (M)$ we will denote by $f^c$ and
$f^v$ the \textit{complete and vertical lift} to $E$ of $f$
defined by
\begin{equation*}
f^c(u)=\sigma (u)(f),\quad f^v(u)=f(\pi (u)),\quad u\in E.
\end{equation*}
For $s\in \Gamma (E)$ we can consider the \textit{vertical lift} of $s$
given by $s^v(u)=s(\pi (u))_u^v,$ for $u\in E,$ where $_u^v:E_{\pi
(u)}\rightarrow T_u(E_{\pi (u)})$ is the canonical isomorphism. There exists
a unique vector field $s^c$ on $E$, the \textit{complete lift} of $s$
satisfying the following conditions:

i) $s^c$ is $\pi $-projectable on $\sigma (s),$

ii) $s^c(\overset{\wedge }{\alpha })=\widehat{\mathcal{L}_s\alpha
},$\\for all $\alpha \in \Gamma (E^{*}),$ where $\overset{\wedge
}{\alpha }(u)=\alpha (\pi (u))(u)$, $u\in E$ (see \cite{Gu1,
Gu2}). \\Considering the prolongation $\mathcal{T}E$ of $E$
\cite{Ma2}, we may introduce the \textit{vertical lift
}$s^{\mathrm{v}}$ and the \textit{complete lift} $s^{ \mathrm{c}}$
of a section $s\in \Gamma (E)$ as the sections of $\mathcal{T}
E\rightarrow E$ given by
\begin{equation*}
s^{\mathrm{v}}(u)=(0,s^v(u)),\quad s^{\mathrm{c}}(u)=(s(\pi
(u)),s^c(u)),\quad u\in E.
\end{equation*}
Other two canonical objects on $\mathcal{T}E$ are the
\textit{Euler section} $\Bbb{C}$ and the \textit{tangent
structure} (\textit{vertical endomorphism}) $J$. The Euler section
$\Bbb{C}$ is the section of $\mathcal{T} E\rightarrow E$ defined
by $\Bbb{C}(u)=(0,u_u^v),\ \forall u\in E.$ The vertical
endomorphism is the section of $(\mathcal{T}E)\oplus (\mathcal{T}
E)^{*}\rightarrow E$ characterized by $J(s^{\mathrm{v}})=0,$
$J(s^{\mathrm{c} })=s^{\mathrm{v}}$, $s\in \Gamma (E)$ which
satisfies
\begin{equation*}
J^2=0,\ ImJ=\ker J=V\mathcal{T}E,\quad \ [\Bbb{C},J]_{\mathcal{T}E}=-J.
\end{equation*}
A section $\mathcal{S}$ of $\mathcal{T}E\rightarrow E$ is called
\textit{ semispray} (\textit{second order differential equation
-SODE}) on $E$ if $J( \mathcal{S})=\Bbb{C}$. The local basis of
$\Gamma (\mathcal{T}E)$ is given by $\{\mathcal{X}_\alpha
,\mathcal{V}_\alpha \}$, where
\begin{equation}
\mathcal{X}_\alpha (u)=\left( s_\alpha (\pi (u)),\left. \sigma
_\alpha ^i\frac \partial {\partial x^i}\right| _u\right) ,\quad
\mathcal{V}_\alpha (u)=\left( 0,\left. \frac \partial {\partial
y^\alpha }\right| _u\right),
\end{equation}
and $(\partial /\partial x^i,\partial /\partial y^\alpha )$ is the
local basis on $TE$ (see \cite{Ma2}). The structure functions of
$\mathcal{T}E$ are given by the following formulas
\begin{equation}
\sigma ^1(\mathcal{X}_\alpha )=\sigma _\alpha ^i\frac \partial {\partial
x^i},\quad \sigma ^1(\mathcal{V}_\alpha )=\frac \partial {\partial y^\alpha
},
\end{equation}
\begin{equation}
\lbrack \mathcal{X}_\alpha ,\mathcal{X}_\beta
]_{\mathcal{T}E}=L_{\alpha \beta }^\gamma \mathcal{X}_\gamma
,\quad [\mathcal{X}_\alpha ,\mathcal{V} _\beta
]_{\mathcal{T}E}=0,\quad [\mathcal{V}_\alpha ,\mathcal{V}_\beta
]_{ \mathcal{T}E}=0.
\end{equation}
The vertical lift of a section $\rho =\rho ^\alpha s_\alpha $ is
$\rho ^{ \mathrm{v}}=\rho ^\alpha \mathcal{V}_\alpha $.$\,$ The
coordinate expression of Euler section is $\Bbb{C}=y^\alpha
\mathcal{V}_\alpha $ and the local expression of $J$ is given by
$J=\mathcal{X}^\alpha \otimes \mathcal{V} _\alpha ,$ where
$\{\mathcal{X}^\alpha ,\mathcal{V}^\alpha \}$ denotes the
corresponding dual basis of $\{\mathcal{X}_\alpha
,\mathcal{V}_\alpha \}$. The Nijenhuis tensor of the vertical
endomorphism vanishes and it results that $J$ is integrable. The
expression of the complete lift of a section $ \rho =\rho ^\alpha
s_\alpha $ is
\begin{equation}
\rho ^{\mathrm{c}}=\rho ^\alpha \mathcal{X}_\alpha +(\sigma
_\varepsilon ^i \frac{\partial \rho ^\alpha }{\partial
x^i}-L_{\beta \varepsilon }^\alpha \rho ^\beta )y^\varepsilon
\mathcal{V}_\alpha .
\end{equation}
In particular $s_\alpha ^{\mathrm{v}}=\mathcal{V}_\alpha $,
$s_\alpha ^{ \mathrm{c}}=\mathcal{X}_\alpha -L_{\alpha \varepsilon
}^\beta y^\varepsilon \mathcal{V}_\beta .$ The local expression of
the differential of a function $ L$ on $\mathcal{T}E$ is
$d^EL=\sigma _\alpha ^i\frac{\partial L}{\partial x^i
}\mathcal{X}^\alpha +\frac{\partial L}{\partial y^\alpha
}\mathcal{V}^\alpha $ and we have $d^Ex^i=\sigma _\alpha
^i\mathcal{X}^\alpha $, $\ d^Ey^\alpha = \mathcal{V}^\alpha$. The
differential of sections of $(\mathcal{T}E)^{*}$ is determined by
\begin{equation*}
d^E\mathcal{X}^\alpha =-\frac 12L_{\beta \gamma }^\alpha \mathcal{X}^\beta
\wedge \mathcal{X}^\gamma ,\quad d^E\mathcal{V}^\alpha =0.
\end{equation*}
In local coordinates a semispray has the expression
\begin{equation}
\mathcal{S}(x,y)=y^\alpha \mathcal{X}_\alpha +\mathcal{S}^\alpha
(x,y) \mathcal{V}_\alpha .
\end{equation}
and the following equality holds
\begin{equation}
J[\mathcal{S},JX]_{\mathcal{T}E}=-JX,\ X\in \Gamma (E).
\end{equation}
The integral curves of $\sigma ^1(\mathcal{S})$ satisfy the differential
equations
\begin{equation*}
\frac{dx^i}{dt}=\sigma _\alpha ^i(x)y^\alpha ,\quad
\frac{dy^\alpha }{dt}= \mathcal{S}^\alpha (x,y).
\end{equation*}
If we have the relation
$[\Bbb{C},\mathcal{S}]_{\mathcal{T}E}=\mathcal{S}$ then
$\mathcal{S}$ is called spray and the functions
$\mathcal{S}^\alpha $ are homogeneous functions of degree $2$ in
$y^\alpha .$  Let us consider a regular Lagrangian $L$ on $E$,
that is the matrix
\begin{equation*}
g_{\alpha \beta }=\frac{\partial ^2L}{\partial y^\alpha \partial
y^\beta},
\end{equation*}
has constant rank $m$. We have the Cartan 1-section $\theta
_L=\frac{\partial L}{\partial y^\alpha }\mathcal{X}^\alpha$ and
the Cartan 2-section $\omega _L=d^E\theta _L$, which is a
symplectic structure induced by $L$ given by \cite{Ma2}
\begin{equation*}
\omega _L=g_{\alpha \beta }\mathcal{V}^\beta \wedge
\mathcal{X}^\alpha +\frac 12\left( \sigma _\alpha ^i\frac{\partial
^2L}{\partial x^i\partial y^\beta }-\sigma _\beta ^i\frac{\partial
^2L}{\partial x^i\partial y^\alpha } -\frac{\partial L}{\partial
y^\varepsilon }L_{\alpha \beta }^\varepsilon \right)
\mathcal{X}^\alpha \wedge \mathcal{X}^\beta .
\end{equation*}
Considering the energy function $E_L=\Bbb{C}(L)-L$, with local expression
\begin{equation*}
E_L=y^\alpha \frac{\partial L}{\partial y^\alpha }-L,
\end{equation*}
then the symplectic equation
\begin{equation*}
i_S\omega _L=-d^EE_L,
\end{equation*}
determines the components of the canonical semispray \cite{Ma2}
\begin{equation}
S^\varepsilon =g^{\varepsilon \beta }\left( \sigma _\beta
^i\frac{\partial L }{\partial x^i}-\sigma _\alpha ^i\frac{\partial
^2L}{\partial x^i\partial y^\beta }y^\alpha -L_{\beta \alpha
}^\gamma y^\alpha \frac{\partial L}{
\partial y^\gamma }\right) ,
\end{equation}
where $g_{\alpha \beta }g^{\beta \gamma }=\delta _\alpha ^\gamma $, which
depends only on the regular Lagrangian and the structure function of the Lie
algebroid.

\subsection{\textbf{Nonlinear connections on Lie algebroids}}

A nonlinear connection is an important tool in the geometry of
systems of second order differential equations. The system of SODE
can by represented using the notion of semispray, which together
with a nonlinear connection induce two important concepts (the
dynamical covariant derivative and Jacobi endomorphism) which are
used in order to find the invariant equations of the symmetries of
SODE.
\begin{definition}
A nonlinear connection on $\mathcal{T}E$ is an almost product
structure $ \mathcal{N}$ on $\pi _2:\mathcal{T}E\rightarrow E$
(i.e. a bundle morphism $ \mathcal{N}:\mathcal{T}E\rightarrow
\mathcal{T}E$, such that $\mathcal{N} ^2=Id$) smooth on
$\mathcal{T}E\backslash \{0\}$ such that $V\mathcal{T} E=\ker
(Id+\mathcal{N}).$
\end{definition}
If $\mathcal{N}$ is a connection on $\mathcal{T}E$ then
$H\mathcal{T}E=\ker (Id-\mathcal{N})$ is the horizontal subbundle
associated to $\mathcal{N}$ and $\mathcal{T}E=V\mathcal{T}E\oplus
H\mathcal{T}E.$ Each $\rho \in \Gamma ( \mathcal{T}E)$ can be
written as $\rho =\rho ^{\mathrm{h}}+\rho ^{\mathrm{v} },$ where
$\rho ^{\mathrm{h}}$, $\rho ^{\mathrm{v}}$ are sections in the
horizontal and respective vertical subbundles. If $\rho
^{\mathrm{h}}=0,$ then $\rho $ is called\textit{\ vertical }and if
$\rho ^{\mathrm{v}}=0,$ then $\rho $ is called
\textit{horizontal}. A connection $\mathcal{N}$ on $ \mathcal{T}E$
induces two projectors $\mathrm{h},\mathrm{v}:\mathcal{T}
E\rightarrow \mathcal{T}E$ such that $\mathrm{h}(\rho )=\rho
^{\mathrm{h}}$ and $\mathrm{v}(\rho )=\rho ^{\mathrm{v}}$ for
every $\rho \in \Gamma ( \mathcal{T}E)$. We have
\begin{equation*}
\mathrm{h}=\frac 12(Id+\mathcal{N}),\quad \mathrm{v}=\frac
12(Id-\mathcal{N} ),\quad \ker
\mathrm{h}=Im\mathrm{v}=V\mathcal{T}E,\quad Im\mathrm{h}=\ker
\mathrm{v}=H\mathcal{T}E.
\end{equation*}
\begin{equation*}
\mathrm{h}^2=\mathrm{h},\quad \mathrm{v}^2=\mathrm{v},\quad
\mathrm{hv}= \mathrm{vh}=0,\quad \mathrm{h}+\mathrm{v}=Id,\quad
\mathrm{h}-\mathrm{v}= \mathcal{N}.
\end{equation*}
\begin{equation*}
J\mathrm{h}=J,\quad \mathrm{h}J=0,\quad J\mathrm{v}=0,\quad \mathrm{v}J=J.
\end{equation*}
Locally, a connection can be expressed as
$\mathcal{N}(\mathcal{X}_\alpha )= \mathcal{X}_\alpha
-2\mathcal{N}_\alpha ^\beta \mathcal{V}_\beta $, $
\mathcal{N}(\mathcal{V}_\beta )=-\mathcal{V}_\beta ,$ where
$\mathcal{N} _\alpha ^\beta =\mathcal{N}_\alpha ^\beta (x,y)$ are
\ the local coefficients of $\mathcal{N}$. The sections
\begin{equation*}
\delta _\alpha =\mathrm{h}(\mathcal{X}_\alpha )=\mathcal{X}_\alpha
-\mathcal{ N}_\alpha ^\beta \mathcal{V}_\beta ,
\end{equation*}
generate a basis of $H\mathcal{T}E$. The frame $\{\delta _\alpha
,\mathcal{V} _\alpha \}$ is a local basis of $\mathcal{T}E$ called
Berwald basis. The dual adapted basis is $\{\mathcal{X}^\alpha
,\delta \mathcal{V}^\alpha \}$ where $\delta \mathcal{V}^\alpha
=\mathcal{V}^\alpha -\mathcal{N}_\beta ^\alpha \mathcal{X}^\beta
.$ The Lie brackets of the adapted basis $\{\delta _\alpha
,\mathcal{V}_\alpha \}$ are \cite{Po2}
\begin{equation}
\lbrack \delta _\alpha ,\delta _\beta ]_{\mathcal{T}E}=L_{\alpha
\beta }^\gamma \delta _\gamma +\mathcal{R}_{\alpha \beta }^\gamma
\mathcal{V} _\gamma ,\quad [\delta _\alpha ,\mathcal{V}_\beta
]_{\mathcal{T}E}=\frac{
\partial \mathcal{N}_\alpha ^\gamma }{\partial y^\beta }\mathcal{V}_\gamma
,\quad [\mathcal{V}_\alpha ,\mathcal{V}_\beta ]_{\mathcal{T}E}=0,
\end{equation}
\begin{equation}
\mathcal{R}_{\alpha \beta }^\gamma =\delta _\beta (\mathcal{N}_\alpha
^\gamma )-\delta _\alpha (\mathcal{N}_\beta ^\gamma )+L_{\alpha \beta
}^\varepsilon \mathcal{N}_\varepsilon ^\gamma .
\end{equation}

\begin{definition}
The curvature of the nonlinear connection $\mathcal{N}$ on
$\mathcal{T}E$ is $\Omega =-\mathbf{N}_{\mathrm{h}}$ where
$\mathrm{h}$ is the horizontal projector and
$\mathbf{N}_{\mathrm{h}}$ is the Nijenhuis tensor of $\mathrm{h
}$.
\end{definition}
In local coordinates we have
\begin{equation*}
\Omega =-\frac 12\mathcal{R}_{\alpha \beta }^\gamma \mathcal{X}^\alpha
\wedge \mathcal{X}^\beta \otimes \mathcal{V}_\gamma ,
\end{equation*}
where $\mathcal{R}_{\alpha \beta }^\gamma $ are given by (11) and
represent the local coordinate functions of the curvature tensor.
The curvature of the nonlinear connection is an obstruction to the
integrability of $H\mathcal{T} E $, understanding that a vanishing
curvature entails that horizontal sections are closed under the
Lie algebroid bracket of $\mathcal{T}E$. The horizontal
distribution {}$H\mathcal{T}E$ is integrable if and only if the
curvature $\Omega $ of the nonlinear connection vanishes. Also,
from the Jacobi identity we obtain
\begin{equation*}
\lbrack \mathrm{h},\Omega ]_{\mathcal{T}E}=0.
\end{equation*}
Let us consider a semispray $\mathcal{S}$ and an arbitrary
nonlinear connection $\mathcal{N}$ with induced
$(\mathrm{h},\mathrm{v})$ projectors. Then we set (see also
\cite{Po3}).

\begin{definition}
The vertically valued $(1,1)$-type tensor field on Lie algebroid
$\mathcal{T} E$ given by
\begin{equation}
\Phi =-\mathrm{v}\circ \mathcal{L}_{\mathcal{S}}\mathrm{v},
\end{equation}
will be called the Jacobi endomorphism.
\end{definition}
The Jacobi endomorphism $\Phi$ has been used in the study of
Jacobi equations for SODE on Lie algebroids in \cite{Ca1} and to
express one of the Helmholtz conditions of the inverse problem of
the calculus of variations on Lie algebroids \cite{Po3} (see also
\cite{Bar}).
 We obtain
\begin{equation*}
\Phi =-\mathrm{v}\circ
\mathcal{L}_{\mathcal{S}}\mathrm{v}=\mathrm{v}\circ
\mathcal{L}_{\mathcal{S}}\mathrm{h}=\mathrm{v}\circ
(\mathcal{L}_{\mathcal{S} }\circ \mathrm{h}-\mathrm{h}\circ
\mathcal{L}_{\mathcal{S}})=\mathrm{v}\circ
\mathcal{L}_{\mathcal{S}}\circ \mathrm{h},
\end{equation*}
and in local coordinates the action of Lie derivative on the Berwald basis
is given by
\begin{equation}
\mathcal{L}_{\mathcal{S}}\delta _\beta =\left( \mathcal{N}_\beta
^\alpha -L_{\beta \varepsilon }^\alpha y^\varepsilon \right)
\delta _\alpha + \mathcal{R}_\beta ^\gamma \mathcal{V}_\gamma
,\quad \mathcal{L}_{\mathcal{S}} \mathcal{V}_\beta =-\delta _\beta
-\left( \mathcal{N}_\beta ^\alpha +\frac{
\partial \mathcal{S}^\alpha }{\partial y^\beta }\right) \mathcal{V}_\alpha.
\end{equation}
The Jacobi endomorphism has the local form
\begin{equation}
\Phi =\mathcal{R}_\beta ^\alpha \mathcal{V}_\alpha \otimes
\mathcal{X}^\beta ,\quad \mathcal{R}_\beta ^\gamma =-\sigma _\beta
^i\frac{\partial \mathcal{S} ^\gamma }{\partial
x^i}-\mathcal{S}(\mathcal{N}_\beta ^\gamma )+\mathcal{N} _\beta
^\alpha \mathcal{N}_\alpha ^\gamma +\mathcal{N}_\beta ^\alpha
\frac{
\partial \mathcal{S}^\gamma }{\partial y^\alpha }+\mathcal{N}_\varepsilon
^\gamma L_{\alpha \beta }^\varepsilon y^\alpha .
\end{equation}

\begin{proposition}
The following formula holds
\begin{equation}
\Phi =i_{\mathcal{S}}\Omega +\mathrm{v}\circ
\mathcal{L}_{\mathrm{v}\mathcal{ S}}\mathrm{h}.
\end{equation}
\end{proposition}

\textbf{Proof}. Indeed, $\Phi (\rho )=\mathrm{v}\circ
\mathcal{L}_{\mathcal{S }}\mathrm{h}\rho =\mathrm{v}\circ
\mathcal{L}_{\mathrm{h}\mathcal{S}}\mathrm{ h}\rho
+\mathrm{v}\circ \mathcal{L}_{\mathrm{v}\mathcal{S}}\mathrm{h}\rho
$ and $\Omega (\mathcal{S},\rho
)=\mathrm{v}[\mathrm{h}\mathcal{S},\mathrm{h} \rho
]_{\mathcal{T}E}=\mathrm{v}\circ
\mathcal{L}_{\mathrm{h}\mathcal{S}} \mathrm{h}\rho$, which yields
$\Phi (\rho )=\Omega (\mathcal{S},\rho )+ \mathrm{v}\circ
\mathcal{L}_{\mathrm{v}\mathcal{S}}\mathrm{h}\rho .$ \hfill
\hbox{\rlap{$\sqcap$}$\sqcup$}

For a given semispray $\mathcal{S}$ on $\mathcal{T}E$ the Lie
derivative $\mathcal{L}_\mathcal{S}$ defines a tensor derivation
on $\mathcal{T}E$, but does not preserve some of the geometric
structure as tangent structure and nonlinear connection. Next,
using a nonlinear connection, we introduce a tensor derivation on
$\mathcal{T}E$, called the dynamical covariant derivative, that
preserves some other geometric structures.

\section{\textbf{Dynamical covariant derivative on Lie algebroids}}

In the following we will introduce the notion of dynamical
covariant derivative on Lie algebroids as a tensor derivation and
study its properties. We will use the Jacobi endomorphism and the
dynamical covariant derivative in the study of symmetries for SODE
on Lie algebroids.

\begin{definition}
\cite{Po3} A map $\nabla :\frak{T}(\mathcal{T}E\backslash \{0\})\rightarrow
\frak{T}(\mathcal{T}E\backslash \{0\})$ is said to be a tensor derivation on
$\mathcal{T}E\backslash \{0\}$ if the following conditions are satisfied:\\
i) $\nabla $ is $\Bbb{R}$-linear\\ii) $\nabla $ is type
preserving, i.e. $ \nabla (\frak{T}_s^r(\mathcal{T}E\backslash
\{0\})\subset \frak{T}_s^r(
\mathcal{T}E\backslash \{0\})$, for each $(r,s)\in \Bbb{N}\times \Bbb{N.}$\\
iii) $\nabla $ obeys the Leibnitz rule $\nabla (P\otimes S)=\nabla
P\otimes S+P\otimes \nabla S$, for any tensors $P,S$ on
$\mathcal{T}E\backslash \{0\}. $\\iv) $\nabla \,$commutes with any
contractions, where $\frak{T}_{\bullet }^{\bullet
}(\mathcal{T}E\backslash \{0\})$ is the space of tensors on
$\mathcal{T}E\backslash \{0\}.$
\end{definition}

For a semispray $\mathcal{S}$ and an arbitrary nonlinear
connection $ \mathcal{N}$ we consider the $\Bbb{R}$-linear map
$\nabla :\Gamma (\mathcal{T }E\backslash \{0\})\rightarrow \Gamma
(\mathcal{T}E\backslash \{0\})$ given by
\begin{equation}
\nabla =\mathrm{h}\circ \mathcal{L}_{\mathcal{S}}\circ
\mathrm{h}+\mathrm{v} \circ \mathcal{L}_{\mathcal{S}}\circ
\mathrm{v,}
\end{equation}
which will be called the dynamical covariant derivative induced by
the semispray $\mathcal{S}$ and the nonlinear connection
$\mathcal{N}$. By setting $\nabla f=\mathcal{S}(f),$ for $f\in
C^\infty (E\backslash \{0\})$ using the Leibnitz rule and the
requirement that $\nabla $ commutes with any contraction, we can
extend the action of $\nabla $ to arbitrary section on $
\mathcal{T}E\backslash \{0\}$. For a section on
$\mathcal{T}E\backslash \{0\} $ the dynamical covariant derivative
is given by $(\nabla \varphi )(\rho )=S(\varphi )(\rho )-\varphi
(\nabla \rho ).$ For a $(1,1)$-type tensor field $T$ on
$\mathcal{T}E\backslash \{0\}$ the dynamical covariant derivative
has the form
\begin{equation}
\nabla T=\nabla \circ T-T\circ \nabla.
\end{equation}
and by direct computation using (17) we obtain
\begin{equation*}
\nabla \mathrm{h}=\nabla \mathrm{v}=0.
\end{equation*}
which means that $\nabla$ preserves the horizontal and vertical
sections. Also, we get
\begin{equation*}
\nabla \mathcal{V}_\beta =\mathrm{v}[\mathcal{S},\mathcal{V}_\beta
]_{ \mathcal{T}E}=-\left( \mathcal{N}_\beta ^\alpha
+\frac{\partial \mathcal{S} ^\alpha }{\partial y^\beta }\right)
\mathcal{V}_\alpha ,\quad \nabla \delta \mathcal{V}^\beta =\left(
\mathcal{N}_\alpha ^\beta +\frac{\partial \mathcal{ S}^\beta
}{\partial y^\alpha }\right) \delta \mathcal{V}^\alpha,
\end{equation*}
\begin{equation*}
\nabla \delta _\beta =\mathrm{h}[\mathcal{S},\delta _\beta
]_{\mathcal{T} E}=\left( \mathcal{N}_\beta ^\alpha -L_{\beta
\varepsilon }^\alpha y^\varepsilon \right) \delta _\alpha,\quad
\nabla \mathcal{X}^\beta =-\left( \mathcal{N}_\alpha ^\beta
-L_{\alpha \varepsilon }^\beta y^\varepsilon \right)
\mathcal{X}^\alpha .
\end{equation*}
The action of the dynamical covariant derivative on the horizontal
section $ X=hX$ is given by following relations
\begin{equation}
\nabla X=\nabla \left( X^\alpha \delta _\alpha \right) =\nabla
X^\alpha \delta _\alpha ,\quad \nabla X^\alpha
=\mathcal{S}(X^\alpha )+\left( \mathcal{N}_\beta ^\alpha
+y^\varepsilon L_{\varepsilon \beta }^\alpha \right) X^\beta.
\end{equation}

\begin{proposition}
The following results hold
\begin{equation}
\mathrm{h}\circ \mathcal{L}_{\mathcal{S}}\circ J=-\mathrm{h},\quad J\circ
\mathcal{L}_{\mathcal{S}}\circ \mathrm{v}=-\mathrm{v},
\end{equation}
\begin{equation}
\nabla J=\mathcal{L}_{\mathcal{S}}J+\mathcal{N},\quad \nabla
J=-\left( \frac{
\partial \mathcal{S}^\beta }{\partial y^\alpha }-y^\varepsilon L_{\alpha
\varepsilon }^\beta +2\mathcal{N}_\alpha ^\beta \right)
\mathcal{V}_\beta \otimes \mathcal{X}^\alpha.
\end{equation}
\end{proposition}

\textbf{Proof}. From (8) we get
\begin{equation*}
J[\mathcal{S},JX]_{\mathcal{T}E}=-JX\Rightarrow J\left(
[\mathcal{S},JX]_{ \mathcal{T}E}+X\right) =0\Rightarrow
[\mathcal{S},JX]_{\mathcal{T}E}+X\in V \mathcal{T}E,
\end{equation*}
\begin{equation*}
\mathrm{h}\left( [\mathcal{S},JX]_{\mathcal{T}E}+X\right) =0\Rightarrow
\mathrm{h}[\mathcal{S},JX]_{\mathcal{T}E}=-\mathrm{h}X\Leftrightarrow
\mathrm{h}\circ \mathcal{L}_{\mathcal{S}}\circ J=-\mathrm{h}.
\end{equation*}
Also, in $J[\mathcal{S},JX]_{\mathcal{T}E}+JX=0$ considering
$JX=\mathrm{v}Z$ it results
$J[\mathcal{S},\mathrm{v}Z]_{\mathcal{T}E}=-\mathrm{v}
Z\Leftrightarrow J\circ \mathcal{L}_{\mathcal{S}}\circ
\mathrm{v}=-\mathrm{v} .$ Next
\begin{eqnarray*}
\nabla \circ J &=&\mathrm{h}\circ \mathcal{L}_{\mathcal{S}}\circ
\mathrm{h} \circ J+\mathrm{v}\circ \mathcal{L}_{\mathcal{S}}\circ
\mathrm{v}\circ J=
\mathrm{v}\circ \mathcal{L}_{\mathcal{S}}\circ J= \\
\ &=&(Id-\mathrm{h})\circ \mathcal{L}_{\mathcal{S}}\circ
J=\mathcal{L} _S\circ J-\mathrm{h}\circ
\mathcal{L}_{\mathcal{S}}\circ J=\mathcal{L} _S\circ J+\mathrm{h}.
\end{eqnarray*}
But, on the other hand
\begin{equation*}
J\circ \nabla =J\circ \mathcal{L}_{\mathcal{S}}\circ
\mathrm{h}=J\circ \mathcal{L}_{\mathcal{S}}\circ
(Id-\mathrm{v})=J\circ \mathcal{L}_{\mathcal{S }}-J\circ
\mathcal{L}_{\mathcal{S}}\circ \mathrm{v}=J\circ \mathcal{L}_{
\mathcal{S}}+v.
\end{equation*}
and we obtain
\begin{equation*}
\nabla \circ J-J\circ \nabla =\mathcal{L}_S\circ
J+\mathrm{h}-J\circ
\mathcal{L}_{\mathcal{S}}-\mathrm{v}\Rightarrow \nabla
J=\mathcal{L}_{
\mathcal{S}}J+\mathrm{h}-\mathrm{v}=\mathcal{L}_{\mathcal{S}}J+\mathcal{N}.
\end{equation*}
For the last relation, we have

\begin{eqnarray*}
\nabla J &=&\nabla \left( \mathcal{X}^\beta \otimes
\mathcal{V}_\beta \right) =\nabla \mathcal{X}^\beta \otimes
\mathcal{V}_\beta +\mathcal{X}
^\beta \otimes \nabla \mathcal{V}_\beta \\
\ &=&-\left( \mathcal{N}_\alpha ^\beta -L_{\alpha \varepsilon
}^\beta y^\varepsilon \right) \mathcal{X}^\alpha \otimes
\mathcal{V}_\beta +\mathcal{ X}^\beta \otimes \left(
-\mathcal{N}_\beta ^\alpha -\frac{\partial \mathcal{S
}^\alpha }{\partial y^\beta }\right) \mathcal{V}_\alpha \\
\ &=&-\mathcal{N}_\alpha ^\beta \mathcal{X}^\alpha \otimes
\mathcal{V}_\beta +L_{\alpha \varepsilon }^\beta y^\varepsilon
\mathcal{X}^\alpha \otimes \mathcal{V}_\beta -\mathcal{N}_\beta
^\alpha \mathcal{X}^\beta \otimes \mathcal{V}_\alpha
-\frac{\partial \mathcal{S}^\alpha }{\partial y^\beta }
\mathcal{X}^\beta \otimes \mathcal{V}_\alpha \\
\ &=&\left( L_{\alpha \varepsilon }^\beta y^\varepsilon -\frac{\partial
\mathcal{S}^\beta }{\partial y^\alpha }-2\mathcal{N}_\alpha ^\beta \right)
\mathcal{X}^\alpha \otimes \mathcal{V}_\beta .
\end{eqnarray*}

\hfill\hbox{\rlap{$\sqcap$}$\sqcup$}.

The above proposition leads to the following result:

\begin{theorem}
For a semispray $\mathcal{S}$, an arbitrary nonlinear connection
$\mathcal{N} $ and $\nabla$ the dynamical covariant derivative
induced by $\mathcal{S}$ and $ \mathcal{N}$, the following
conditions are equivalent:\\
$i)$ $\nabla J=0,$\\
$ii)$ $\mathcal{L}_SJ+\mathcal{N}=0,$\\
$iii)$ $\mathcal{N}_\alpha ^\beta =\frac 12\left( -\frac{\partial
\mathcal{S} ^\beta }{\partial y^\alpha }+y^\varepsilon L_{\alpha
\varepsilon }^\beta \right)$.
\end{theorem}

\textbf{Proof}. The proof follows from the relations (20).\hfill
\hbox{\rlap{$\sqcap$}$\sqcup$}

This theorem shows that the compatibility condition $\nabla J=0$ of the
dynamical covariant derivative with the tangent structure determines the
nonlinear connection $\mathcal{N}=-\mathcal{L}_{\mathcal{S}}J$. For the
particular case of tangent bundle we obtain the results from \cite{Bu3}. In
the following we deal with this nonlinear connection induced by semispray.

\subsection{The canonical nonlinear connection induced by a semispray}

A semispray $\mathcal{S},$ together with the condition $\nabla
J=0,$ determines the canonical nonlinear connection
$\mathcal{N}=-\mathcal{L}_{ \mathcal{S}}J$ with local coefficients
\begin{equation*}
\mathcal{N}_\alpha ^\beta =\frac 12\left( -\frac{\partial \mathcal{S}^\beta
}{\partial y^\alpha }+y^\varepsilon L_{\alpha \varepsilon }^\beta \right) .
\end{equation*}
In this case the following equations hold
\begin{equation*}
\lbrack \mathcal{S},\mathcal{V}_\beta ]_{\mathcal{T}E}=-\delta _\beta
+\left( \mathcal{N}_\beta ^\alpha -L_{\beta \varepsilon }^\alpha
y^\varepsilon \right) \mathcal{V}_\alpha ,
\end{equation*}
\begin{equation*}
\lbrack \mathcal{S},\delta _\beta ]_{\mathcal{T}E}=\left(
\mathcal{N}_\beta ^\alpha -L_{\beta \varepsilon }^\alpha
y^\varepsilon \right) \delta _\alpha + \mathcal{R}_\beta ^\alpha
\mathcal{V}_\alpha ,
\end{equation*}
where
\begin{equation}
\mathcal{R}_\beta ^\alpha =-\sigma _\beta ^i\frac{\partial
\mathcal{S} ^\alpha }{\partial x^i}-\mathcal{S}(\mathcal{N}_\beta
^\alpha )-\mathcal{N} _\gamma ^\alpha \mathcal{N}_\beta ^\gamma
+(L_{\varepsilon \beta }^\gamma \mathcal{N}_\gamma ^\alpha
+L_{\gamma \varepsilon }^\alpha \mathcal{N}_\beta ^\gamma
)y^\varepsilon .
\end{equation}
are the local coefficients of the Jacobi endomorphism.
\begin{proposition}
If $\mathcal{S}$ is a spray, then the Jacobi endomorphism is the
contraction with $\mathcal{S}$ of curvature of the nonlinear
connection
\begin{equation*}
\Phi =i_{\mathcal{S}}\Omega .
\end{equation*}
\end{proposition}

\textbf{Proof}. If $\mathcal{S}$ is a spray, then the coefficients
$\mathcal{ S}^\alpha $ are 2-homogeneous with respect to the
variables $y^\beta$ and it results
\begin{equation*}
2\mathcal{S}^\alpha =\frac{\partial \mathcal{S}^\alpha }{\partial
y^\beta } y^\beta =-2\mathcal{N}_\beta ^\alpha y^\beta +L_{\beta
\gamma }^\alpha y^\beta y^\gamma =-2\mathcal{N}_\beta ^\alpha
y^\beta .
\end{equation*}
\begin{equation*}
\mathcal{S}=\mathrm{h}\mathcal{S}=y^\alpha \delta _\alpha ,\quad
\mathrm{v} \mathcal{S}=0,\quad \mathcal{N}_\beta ^\alpha
=\frac{\partial \mathcal{N} _\varepsilon ^\alpha }{\partial
y^\beta }y^\varepsilon +L_{\beta \varepsilon }^\alpha
y^\varepsilon ,
\end{equation*}
which together with (15) yields $\Phi =i_{\mathcal{S}}\Omega$.
Locally, we get $\mathcal{R} _\beta ^\alpha
=\mathcal{R}_{\varepsilon \beta }^\alpha y^\varepsilon $ and
represents the local relation between the Jacobi endomorphism and
the curvature of the nonlinear connection. Also, we have
$\Phi(\mathcal{S})=0$. \hfill
\hbox{\rlap{$\sqcap$}$\sqcup$}\\
Next, we introduce the almost complex structure in order to find
the decomposition formula for the dynamical covariant derivative.
\begin{definition}
The almost complex structure is given by the formula
\begin{equation*}
\Bbb{F}=\mathrm{h}\circ \mathcal{L}_{\mathcal{S}}\mathrm{h}-J.
\end{equation*}
\end{definition}

We have to show that $\Bbb{F}^2=-Id$. Indeed, from the relation
$\mathcal{L}
_{\mathcal{S}}\mathrm{h}=\mathcal{L}_{\mathcal{S}}\circ
\mathrm{h}-\mathrm{h} \circ \mathcal{L}_{\mathcal{S}}$ we obtain
$\Bbb{F}=\mathrm{h}\circ \mathcal{L}_{\mathcal{S}}\circ
\mathrm{h}- \mathrm{h}\circ
\mathcal{L}_{\mathcal{S}}-J=\mathrm{h}\circ \mathcal{L}_{
\mathcal{S}}\circ (\mathrm{h}-Id)-J=-\mathrm{h}\circ
\mathcal{L}_{\mathcal{S} }\circ \mathrm{v}-J$ and
$\Bbb{F}^2=\left( -\mathrm{h}\circ \mathcal{L}_{\mathcal{S}}\circ
\mathrm{v} -J\right) \circ \left( -\mathrm{h}\circ
\mathcal{L}_{\mathcal{S}}\circ \mathrm{v}-J\right)
=\mathrm{h}\circ \mathcal{L}_{\mathcal{S}}\circ \mathrm{v }\circ
\mathrm{h}\circ \mathcal{L}_{\mathcal{S}}\circ
\mathrm{v}+\mathrm{h} \circ \mathcal{L}_{\mathcal{S}}\circ
\mathrm{v}\circ J+$\\
$+J\circ \mathrm{h}\circ \mathcal{L}_{\mathcal{S}}\circ
\mathrm{v}+J^2= \mathrm{h}\circ \mathcal{L}_{\mathcal{S}}\circ
J+J\circ \mathcal{L}_{ \mathcal{S}}\circ \mathrm{v}$
$=-\mathrm{h}-\mathrm{v}=-Id.$

\begin{proposition}
The following results hold
\begin{equation*}
\begin{array}{c}
\Bbb{F}\circ J=\mathrm{h},\quad J\circ \Bbb{F}=\mathrm{v},\quad
\mathrm{v}
\circ \Bbb{F}=\Bbb{F}\circ \mathrm{h}=-J, \\
\mathrm{h}\circ \Bbb{F}=\Bbb{F}\circ \mathrm{v}=\Bbb{F}+J,\quad
\mathcal{N} \circ \Bbb{F}=\Bbb{F}+2J,\quad \Phi
=\mathcal{L}_{\mathcal{S}}\mathrm{h}- \Bbb{F}-J.
\end{array}
\end{equation*}
\end{proposition}

\textbf{Proof}. Using the relations (19) we obtain\\
 $\Bbb{F}\circ
J=\left( -\mathrm{h}\circ \mathcal{L}_{\mathcal{S}}\circ
\mathrm{v}-J\right) \circ J=$ $-\mathrm{h}\circ
\mathcal{L}_{\mathcal{S} }\circ \mathrm{v}\circ J-J^2$
$=-\mathrm{h}\circ \mathcal{L}_{\mathcal{S} }\circ
J=\mathrm{h}$,\\
$J\circ \Bbb{F}=-J\circ \left( \mathrm{h}\circ
\mathcal{L}_{\mathcal{S} }\circ \mathrm{v}+J\right) =-J\circ
\mathrm{h}\circ \mathcal{L}_{\mathcal{S} }\circ
\mathrm{v}-J^2=-J\circ \mathcal{L}_{\mathcal{S}}\circ \mathrm{v}=
\mathrm{v}$,\\
$\mathrm{v}\circ \Bbb{F}=\mathrm{v}\circ \left( \mathrm{h}\circ
\mathcal{L}_{ \mathcal{S}}\mathrm{h}-J\right) =-\mathrm{v}\circ
J=-J$, $\Bbb{F}\circ \mathrm{h}=\left( -\mathrm{h}\circ
\mathcal{L}_{\mathcal{S}}\circ \mathrm{v} -J\right) \circ
\mathrm{h}=-J\circ \mathrm{h}=-J,$ $\mathrm{h}\circ
\Bbb{F}=\mathrm{h}\circ \left( \mathrm{h}\circ \mathcal{L}_{
\mathcal{S}}\mathrm{h}-J\right) =\mathrm{h}\circ
\mathcal{L}_{\mathcal{S}} \mathrm{h}=\Bbb{F}+J$, $\Bbb{F}\circ
\mathrm{v=}\left( -\mathrm{h}\circ \mathcal{L}_{\mathcal{S}}\circ
\mathrm{v}-J\right) \circ \mathrm{v}=-\mathrm{ h}\circ
\mathcal{L}_{\mathcal{S}}\circ \mathrm{v}=$ $\Bbb{F}+J$. In the
same way, the other relations can be proved.\hfill
\hbox{\rlap{$\sqcap$}$\sqcup$}\\
In local coordinates we have
\begin{equation*}
\Bbb{F}=-\mathcal{V}_\alpha \otimes \mathcal{X}^\alpha
+\delta_\alpha \otimes \delta \mathcal{V}^\alpha .
\end{equation*}
For a semispray $\mathcal{S}$ and the associated nonlinear connection we
consider the $\Bbb{R}$-linear map $\nabla _0:\Gamma (\mathcal{T}E\backslash
\{0\})\rightarrow \Gamma (\mathcal{T}E\backslash \{0\})$ given by
\begin{equation*}
\nabla _0\rho =\mathrm{h}[\mathcal{S},\mathrm{h}\rho
]_{\mathcal{T}E}+ \mathrm{v}[\mathcal{S},\mathrm{v}\rho
]_{\mathcal{T}E},\quad \forall \rho \in \Gamma
(\mathcal{T}E\backslash \{0\}).
\end{equation*}
It results that
\begin{equation*}
\nabla _0(f\rho )=\mathcal{S}(f)\rho +f\nabla _0\rho ,\quad
\forall f\in C^\infty (E),\ \rho \in \Gamma
(\mathcal{T}E\backslash \{0\}).
\end{equation*}
Any tensor derivation on $\mathcal{T}E\backslash \{0\}$ is
completely determined by its actions on smooth functions and
sections on $\mathcal{T} E\backslash \{0\}$ (see \cite{Sz2}
generalized Willmore's theorem). Therefore, there exists a unique
tensor derivation $\nabla $ on $\mathcal{T} E\backslash \{0\}$
such that
\begin{equation*}
\nabla \mid _{C^\infty (E)}=\mathcal{S},\quad \nabla \mid _{\Gamma
(\mathcal{ T}E\backslash \{0\})}=\nabla _0.
\end{equation*}
We will call the tensor derivation $\nabla $, the
\textit{dynamical covariant derivative} induced by the semispray
$\mathcal{S}$ (see \cite{Bu2} for the tangent bundle case).
\begin{proposition}
The dynamical covariant derivative has the following decomposition
\begin{equation}
\nabla =\mathcal{L}_{\mathcal{S}}+\Bbb{F}+J-\Phi.
\end{equation}
\end{proposition}

\textbf{Proof}. Using the formula (16) and the expressions of
$\Bbb{F}$ and $\Phi $ we obtain
\begin{eqnarray*}
\nabla &=&\mathrm{h}\circ \mathcal{L}_{\mathcal{S}}\circ
\mathrm{h}+\mathrm{v
}\circ \mathcal{L}_{\mathcal{S}}\circ \mathrm{v}= \\
&=&\mathrm{h}\circ \left(
\mathcal{L}_{\mathcal{S}}\mathrm{h}+\mathrm{h} \circ
\mathcal{L}_{\mathcal{S}}\right) +\mathrm{v}\circ \left(
\mathcal{L}_{
\mathcal{S}}\mathrm{v}+\mathrm{v}\circ \mathcal{L}_{\mathcal{S}}\right) = \\
&=&\mathrm{h}\circ
\mathcal{L}_{\mathcal{S}}\mathrm{h}+\mathrm{v}\circ
\mathcal{L}_{\mathcal{S}}\mathrm{v}+(\mathrm{h}+\mathrm{v})\circ
\mathcal{L}
_{\mathcal{S}}=\mathcal{L}_{\mathcal{S}}+\mathrm{h}\circ
\mathcal{L}_{ \mathcal{S}}\mathrm{h}+\mathrm{v}\circ
\mathcal{L}_{\mathcal{S}}\mathrm{v}=
\\
&=&\mathcal{L}_{\mathcal{S}}+\Bbb{F}+J-\Phi .
\end{eqnarray*}
In this case the dynamical covariant derivative is characterized
by the following formulas
\begin{equation*}
\nabla \mathcal{V}_\beta =\mathrm{v}[\mathcal{S},\mathcal{V}_\beta
]_{ \mathcal{T}E}=\left( \mathcal{N}_\beta ^\alpha -L_{\beta
\varepsilon }^\alpha y^\varepsilon \right) \mathcal{V}_\alpha
=-\frac 12\left( \frac{
\partial \mathcal{S}^\alpha }{\partial y^\beta }+L_{\beta \varepsilon
}^\alpha y^\varepsilon \right) \mathcal{V}_\alpha ,
\end{equation*}
\begin{equation*}
\nabla \delta _\beta =\mathrm{h}[\mathcal{S},\delta _\beta
]_{\mathcal{T} E}=\left( \mathcal{N}_\beta ^\alpha -L_{\beta
\varepsilon }^\alpha y^\varepsilon \right) \delta _\alpha =-\frac
12\left( \frac{\partial \mathcal{S}^\alpha }{\partial y^\beta
}+L_{\beta \varepsilon }^\alpha y^\varepsilon \right) \delta
_\alpha.
\end{equation*}
The next result shows that $\nabla$ acts identically on both
vertical and horizontal distribution, that is enough to find the
action of $\nabla$ on either one of the two distributions.
\begin{proposition}
The dynamical covariant derivative induced by the semispray $\mathcal{S}$ is
compatible with $J$ and $\Bbb{F}$, that is
\begin{equation*}
\nabla J=0,\ \nabla \Bbb{F}=0.
\end{equation*}
\end{proposition}

\textbf{Proof}. $\nabla J=0$ follows from (20). Using the formula
$\Bbb{F}=- \mathrm{h}\circ \mathcal{L}_{\mathcal{S}}\circ
\mathrm{v}-J$ and $\nabla \Bbb{F}=\nabla \circ
\Bbb{F}-\Bbb{F}\circ \nabla$ we obtain
\begin{eqnarray*}
\nabla \Bbb{F} &=&(\mathrm{h}\circ \mathcal{L}_{\mathcal{S}}\circ
\mathrm{h}+ \mathrm{v}\circ \mathcal{L}_{\mathcal{S}}\circ
\mathrm{v})\circ (-\mathrm{h} \circ \mathcal{L}_{\mathcal{S}}\circ
\mathrm{v})-(-\mathrm{h}\circ \mathcal{L }_{\mathcal{S}}\circ
\mathrm{v})\circ (\mathrm{h}\circ \mathcal{L}_{\mathcal{ S}}\circ
\mathrm{h}+\mathrm{v}\circ \mathcal{L}_{\mathcal{S}}\circ
\mathrm{v}
)= \\
\ &=&-\mathrm{h}\circ \mathcal{L}_{\mathcal{S}}\circ
\mathrm{h}\circ \mathcal{L}_{\mathcal{S}}\circ
\mathrm{v}+\mathrm{h}\circ \mathcal{L}_{ \mathcal{S}}\circ
\mathrm{v}\circ \mathcal{L}_{\mathcal{S}}\circ \mathrm{v}=
\\
\ &=&\mathrm{h}\circ \mathcal{L}_{\mathcal{S}}\circ
(\mathrm{v}-\mathrm{h} )\circ \mathcal{L}_{\mathcal{S}}\circ
\mathrm{v}=\mathrm{h}\circ \mathcal{L} _{\mathcal{S}}\circ
\mathcal{L}_{\mathcal{S}}J\circ \mathcal{L}_{\mathcal{S}
}\circ \mathrm{v}= \\
\ &=&\mathrm{h}\circ \mathcal{L}_{\mathcal{S}}\circ
(\mathcal{L}_{\mathcal{S} }\circ J-J\circ
\mathcal{L}_{\mathcal{S}})\circ \mathcal{L}_{\mathcal{S}
}\circ \mathrm{v}= \\
\ &=&\mathrm{h}\circ \mathcal{L}_{\mathcal{S}}\circ
\mathcal{L}_{\mathcal{S} }\circ (J\circ
\mathcal{L}_{\mathcal{S}}\circ \mathrm{v})-(\mathrm{h}\circ
\mathcal{L}_{\mathcal{S}}\circ J)\circ
\mathcal{L}_{\mathcal{S}}\circ
\mathcal{L}_{\mathcal{S}}\circ \mathrm{v}= \\
\ &=&-\mathrm{h}\circ \mathcal{L}_{\mathcal{S}}\circ
\mathcal{L}_{\mathcal{S} }\circ \mathrm{v}+\mathrm{h}\circ
\mathcal{L}_{\mathcal{S}}\circ \mathcal{L} _{\mathcal{S}}\circ
\mathrm{v}=0.
\end{eqnarray*}
\hfill \hbox{\rlap{$\sqcap$}$\sqcup$}\\
The next proposition proves that in the case of spray $\nabla$ has
more properties.
\begin{proposition}
If the dynamical covariant derivative is induced by a spray $\mathcal{S}$
then
\begin{equation*}
\nabla \mathcal{S}=0,\ \nabla \Bbb{C}=0.
\end{equation*}
\end{proposition}

\textbf{Proof}. Indeed, if $\mathcal{S}$ is a spray then we have
$\mathcal{S}= \mathrm{h}\mathcal{S}$ and $\mathrm{v}\mathcal{S}=0$
and it results $\nabla \mathcal{S}=\mathrm{h}\circ
\mathcal{L}_{\mathcal{S}}\circ \mathrm{h}
\mathcal{S}+\mathrm{v}\circ \mathcal{L}_{\mathcal{S}}\circ
\mathrm{v} \mathcal{S}=\mathrm{h}\circ
\mathcal{L}_{\mathcal{S}}\circ \mathcal{S}=0.$ Also $\nabla
\Bbb{C}=0$ follows from $\mathrm{h}\Bbb{C}=0$, $\mathrm{v}\Bbb{C
}=\Bbb{C}$ and
$[\Bbb{C},\mathcal{S}]_{\mathcal{T}E}=\mathcal{S}$.\hfill
\hbox{\rlap{$\sqcap$}$\sqcup$}\\
Next, we introduce the Berwald linear connection induced by a
nonlinear connection and prove that in the case of homogeneous
second order differential equations it coincides with the
dynamical covariant derivative. The Berwald linear connection is
given by
\begin{equation*}
\mathcal{D}:\Gamma (\mathcal{T}E\backslash \{0\})\times \Gamma
(\mathcal{T} E\backslash \{0\})\rightarrow \Gamma
(\mathcal{T}E\backslash \{0\})
\end{equation*}
\begin{equation*}
\mathcal{D}_XY=\mathrm{v}[\mathrm{h}X,\mathrm{v}Y]_{\mathcal{T}E}+\mathrm{h}[
\mathrm{v}X,\mathrm{h}Y]_{\mathcal{T}E}+J[\mathrm{v}X,(\Bbb{F}+J)Y]_{
\mathcal{T}E}+(\Bbb{F}+J)[\mathrm{h}X,JY]_{\mathcal{T}E}.
\end{equation*}

\begin{proposition}
The Berwald linear connection has the following properties
\begin{equation*}
\mathcal{D}\mathrm{h}=0,\quad \mathcal{D}\mathrm{v}=0,\quad
\mathcal{D} J=0,\quad \mathcal{D}\Bbb{F}=0.
\end{equation*}
\end{proposition}

\textbf{Proof}. Using the properties of the vertical and
horizontal projectors we obtain\\
$\mathcal{D}_X\mathrm{v}Y=\mathrm{v}[\mathrm{h}X,\mathrm{v}Y]_{\mathcal{T}
E}+J[\mathrm{v}X,(\Bbb{F}+J)Y]_{\mathcal{T}E}$ and\\
$\mathrm{v}(\mathcal{D}_XY)=\mathrm{v}[\mathrm{h}X,\mathrm{v}Y]_{\mathcal{T}
E}+J[\mathrm{v}X,(\Bbb{F}+J)Y]_{\mathcal{T}E}$ which yields
$\mathcal{D} \mathrm{v}=0$. Also,\\
$\mathcal{D}_X\mathrm{h}Y=\mathrm{h}[\mathrm{v}X,\mathrm{h}Y]_{\mathcal{T}
E}+(\Bbb{F}+J)[\mathrm{h}X,JY]_{\mathcal{T}E}=\mathrm{h}(\mathcal{D}_XY)$
and it results $\mathcal{D}\mathrm{h}=0$. Moreover,\\
$\mathcal{D}_XJY=\mathrm{v}[\mathrm{h}X,JY]_{\mathcal{T}E}+J[\mathrm{v}X,
\mathrm{h}Y]_{\mathcal{T}E}$ and
$J(\mathcal{D}_XY)=J[\mathrm{v}X,\mathrm{h}
Y]_{\mathcal{T}E}+\mathrm{v}[\mathrm{h}X,JY]_{\mathcal{T}E}$ and
we obtain $\mathcal{D}J=0.$ From\\
$\mathcal{D}_X\Bbb{F}Y=\mathrm{v}[\mathrm{h}X,-JY]_{\mathcal{T}E}+\mathrm{h}[
\mathrm{v}X,(\Bbb{F}+J)Y]_{\mathcal{T}E}+J[\mathrm{v}X,-\mathrm{h}Y]_{
\mathcal{T}E}+(\Bbb{F}+J)[\mathrm{h}X,\mathrm{v}Y]_{\mathcal{T}E}$
and\\
$\Bbb{F}(\mathcal{D}_XY)=(\Bbb{F}+J)[\mathrm{h}X,\mathrm{v}Y]_{\mathcal{T}
E}-J[\mathrm{v}X,\mathrm{h}Y]_{\mathcal{T}E}+\mathrm{h}[\mathrm{v}X,(\Bbb{F}
+J)Y]_{\mathcal{T}E}-v[\mathrm{h}X,JY]_{\mathcal{T}E}=$\\
$\mathcal{D}_X\Bbb{F}Y$ we get $\mathcal{D}\Bbb{F}=0.$ \hfill
\hbox{\rlap{$\sqcap$}$\sqcup$}\\
It results that the Berwald connection preserves both horizontal
and vertical sections. Moreover, $\mathcal{D}$ has the same action
on horizontal and vertical distributions and locally we have the
following formulas
\begin{equation*}
\mathcal{D}_{\delta _\alpha }\delta _\beta =\frac{\partial
\mathcal{N} _\alpha ^\gamma }{\partial y^\beta }\delta _\gamma
,\quad \mathcal{D} _{\delta _\alpha }\mathcal{V}_\beta
=\frac{\partial \mathcal{N}_\alpha ^\gamma }{\partial y^\beta
}\mathcal{V}_\gamma ,\quad \mathcal{D}_{\mathcal{V }_\alpha
}\delta _\beta =0,\quad \mathcal{D}_{\mathcal{V}_\alpha
}\mathcal{V} _\beta =0.
\end{equation*}
We can see that the dynamical covariant derivative has the same
properties and this leads to the next result.
\begin{proposition}
If $\mathcal{S}$ is a spray then the following equality holds
\begin{equation*}
\nabla =\mathcal{D}_{\mathcal{S}}.
\end{equation*}
\end{proposition}
\textbf{Proof}. If $\mathcal{S}$ is a spray then
$\mathcal{S}=\mathrm{h} \mathcal{S}$ and $\mathrm{v}\mathcal{S}=0$
which implies
\begin{equation*}
\mathcal{D}_{\mathcal{S}}Y=\mathrm{v}[\mathcal{S},\mathrm{v}Y]_{\mathcal{T}
E}+(\Bbb{F}+J)[\mathcal{S},JY]_{\mathcal{T}E}.
\end{equation*}
But $\nabla
Y=\mathrm{h}[\mathcal{S},\mathrm{h}Y]_{\mathcal{T}E}+\mathrm{v}[
\mathcal{S},\mathrm{v}Y]_{\mathcal{T}E}$ and we will prove that
$\mathrm{h}[
\mathcal{S},\mathrm{h}Y]_{\mathcal{T}E}=(\Bbb{F}+J)[\mathcal{S},JY]_{
\mathcal{T}E}$ using the computation in local coordinates. Let us
consider $ Y=X^\alpha (x,y)\mathcal{X}_\alpha +Y^\beta
(x,y)\mathcal{V}_\beta $ and using (10) we get
\begin{equation*}
\lbrack \mathcal{S},\mathrm{h}Y]_{\mathcal{T}E}=[y^\alpha \delta
_\alpha ,X^\beta \delta _\beta ]_{\mathcal{T}E}=y^\alpha X^\beta
\mathcal{R}_{\alpha \beta }^\varepsilon \mathcal{V}_\varepsilon
+y^\alpha X^\beta L_{\alpha \beta }^\varepsilon \delta
_\varepsilon +y^\alpha \delta _\alpha (X^\beta )\delta _\beta
+X^\beta N_\beta ^\alpha \delta _\alpha,
\end{equation*}

\begin{equation*}
\mathrm{h}[\mathcal{S},\mathrm{h}Y]_{\mathcal{T}E}=\left( y^\alpha
\delta _\alpha (X^\beta )+X^\alpha N_\alpha ^\beta+y^\alpha
X^\varepsilon L_{\alpha \varepsilon}^\beta \right) \delta _\beta.
\end{equation*}
Next
\begin{equation*}
\lbrack \mathcal{S},JY]_{\mathcal{T}E}=[y^\alpha \delta _\alpha ,X^\beta
\mathcal{V}_\beta ]_{\mathcal{T}E}=y^\alpha X^\beta \frac{\partial N_\alpha
^\varepsilon }{\partial y^\beta }\mathcal{V}_\varepsilon +y^\alpha \delta
_\alpha (X^\beta )\mathcal{V}_\beta -X^\beta \delta _\beta .
\end{equation*}
Also, we have
\begin{equation*}
y^\alpha X^\beta \frac{\partial N_\alpha ^\varepsilon }{\partial
y^\beta }= \mathcal{N}_\beta ^\varepsilon X^\beta -L_{\beta \alpha
}^\varepsilon y^\alpha X^\beta,
\end{equation*}
and using the relations $(\Bbb{F}+J)(\mathcal{V}_\alpha )=\delta
_\alpha $, $ (\Bbb{F}+J)(\delta _\alpha )=0$ we obtain the result
which ends the proof. \hfill \hbox{\rlap{$\sqcap$}$\sqcup$}\\
Moreover, $\nabla{\mathcal{S}}
=\mathcal{D}_{\mathcal{S}}{\mathcal{S}}=0$ and it results that the
integral curves of the spray are geodesics of the Berwald linear
connection.

\section{\textbf{Symmetries for semispray}}

In this section we study the symmetries of SODE on Lie algebroids
and prove that the canonical nonlinear connection can be
determined by these symmetries. We find the relations between
dynamical symmetries, Lie symmetries, Newtonoid sections, Cartan
symmetries and conservation laws, and show when one of them will
imply the others. Also, we obtain the invariant equations of these
symmetries, using the dynamical covariant derivative and Jacobi
endomorphism. In the particular case of the tangent bundle some
results from \cite{Bu3, Mar, Pr1, Pr2} are obtained.

\begin{definition}
A section $X\in \Gamma (\mathcal{T}E\backslash \{0\})$ is a dynamical
symmetry of semispray $\mathcal{S}$ if $[\mathcal{S},X]_{\mathcal{T}E}=0.$
\end{definition}

In local coordinates for $X=X^\alpha (x,y)\mathcal{X}_\alpha
+Y^\alpha (x,y) \mathcal{V}_\alpha $ we obtain
\begin{equation*}
\lbrack \mathcal{S},X]_{\mathcal{T}E}=\left( y^\alpha L_{\alpha \gamma
}^\beta X^\gamma -Y^\beta +\mathcal{S}(X^\beta )\right) \mathcal{X}_\beta
+\left( \mathcal{S}(Y^\beta )-X(\mathcal{S}^\beta )\right) \mathcal{V}_\beta
,
\end{equation*}
and it results that the dynamical symmetry is characterized by the equations
\begin{equation}
Y^\alpha =\mathcal{S}(X^\alpha )+y^\varepsilon L_{\varepsilon \beta }^\alpha
X^\beta ,
\end{equation}
\begin{equation}
\mathcal{S}(Y^\alpha )-X(\mathcal{S}^\alpha )=0.
\end{equation}
Introducing (23) into (24) we obtain
\begin{equation*}
\mathcal{S}^2(X^\alpha )-X(\mathcal{S}^\alpha )=\left( \sigma
_\gamma ^i \frac{\partial L_{\varepsilon \beta }^\alpha }{\partial
x^i}X^\beta +L_{\varepsilon \beta }^\alpha \sigma _\gamma
^i\frac{\partial X^\beta }{
\partial x^i}\right) y^\gamma y^\varepsilon +\mathcal{S}^\gamma \left(
L_{\gamma \beta }^\alpha X^\beta +y^\varepsilon L_{\varepsilon
\beta }^\alpha \frac{\partial X^\beta }{\partial y^\gamma
}\right).
\end{equation*}

\begin{definition}
A section $\widetilde{X}=\widetilde{X}^\alpha (x,y)s_\alpha $ on
$E\backslash \{0\}$ is a Lie symmetry of a semispray if its
complete lift $\widetilde{X}^c$ is a dynamical symmetry, that is
$[\mathcal{S},\widetilde{X}^c]_{\mathcal{T}E}=0.$
\end{definition}

\begin{proposition}
The local expression of a Lie symmetry is given by
\begin{eqnarray*}
\mathcal{S}^\alpha \frac{\partial \widetilde{X}^\beta }{\partial
y^\alpha } =0,
\end{eqnarray*}
\begin{equation*}
\mathcal{S}^\alpha \widetilde{X}_{\mid _\alpha }^\beta +y^\alpha
y^\varepsilon \sigma _\alpha ^i\frac{\partial \widetilde{X}_{\mid
_\varepsilon }^\beta }{\partial x^i}-\widetilde{X}^\alpha \sigma
_\alpha ^i \frac{\partial \mathcal{S}^\beta }{\partial
x^i}-y^\varepsilon \widetilde{X} _{\mid _\varepsilon }^\alpha
\frac{\partial \mathcal{S}^\beta }{\partial y^\alpha }+S^\alpha
y^\varepsilon \left( \sigma _\varepsilon ^i\frac{
\partial ^2\widetilde{X}^\beta }{\partial y^\alpha \partial x^i}-L_{\gamma
\varepsilon }^\beta \frac{\partial \widetilde{X}^\gamma }{\partial
y^\alpha } \right) =0.
\end{equation*}
where
\begin{equation*}
\widetilde{X}_{\mid_\varepsilon }^\alpha :=\sigma_\varepsilon
^i\frac{
\partial \widetilde{X}^\alpha }{\partial x^i}-L_{\beta \varepsilon}^\alpha
\widetilde{X}^\beta ,
\end{equation*}
\end{proposition}

\textbf{Proof}. Considering $\widetilde{X}^c=\widetilde{X}^\alpha
\mathcal{X} _\alpha +y^\varepsilon \widetilde{X}_{\mid \varepsilon
}^\alpha \mathcal{V} _\alpha $ and using (1) we obtain

\begin{equation*}
\lbrack \mathcal{S},\widetilde{X}^c]_{\mathcal{T}E}=\left(
\widetilde{X} ^\alpha y^\varepsilon L_{\varepsilon \alpha }^\beta
+y^\alpha \sigma _\alpha ^i\frac{\partial \widetilde{X}^\beta
}{\partial x^i}-y^\varepsilon \widetilde{X}_{\mid _\varepsilon
}^\beta +\mathcal{S}^\alpha \frac{\partial \widetilde{X}^\beta
}{\partial y^\alpha }\right) \mathcal{X}_\beta +\qquad \qquad
\qquad \qquad \qquad \qquad \qquad \qquad \qquad \qquad
\end{equation*}

\begin{equation*}
\ \ \ \ \left( y^\alpha y^\varepsilon \sigma _\alpha
^i\frac{\partial \widetilde{X}_{\mid _\varepsilon }^\beta
}{\partial x^i}-\widetilde{X} ^\alpha \sigma _\alpha
^i\frac{\partial \mathcal{S}^\beta }{\partial x^i}+
\mathcal{S}^\alpha \widetilde{X}_{\mid _\alpha }^\beta
-y^\varepsilon \widetilde{X}_{\mid _\varepsilon }^\alpha
\frac{\partial \mathcal{S}^\beta }{
\partial y^\alpha }+S^\alpha y^\varepsilon \left( \sigma _\varepsilon ^i
\frac{\partial ^2\widetilde{X}^\beta }{\partial y^\alpha \partial
x^i} -L_{\gamma \varepsilon }^\beta \frac{\partial
\widetilde{X}^\gamma }{
\partial y^\alpha }\right) \right) \mathcal{V}_\beta.
\end{equation*}
We deduce that $\widetilde{X}^\alpha y^\varepsilon L_{\varepsilon
\alpha }^\beta +y^\alpha \sigma _\alpha ^i\frac{\partial
\widetilde{X}^\beta }{
\partial x^i}-y^\varepsilon \widetilde{X}_{\mid _\varepsilon }^\beta =0$ and
it results the local expression of a Lie symmetry. \hfill
\hbox{\rlap{$\sqcap$}$\sqcup$}\\
We have to remark that a section $\widetilde{X}=\widetilde{X}^\alpha
(x)s_\alpha $ on $E\backslash \{0\}$ is a Lie symmetry if and only if (see
also \cite{Pe})
\begin{equation*}
y^\alpha y^\varepsilon \sigma _\alpha ^i\frac{\partial
\widetilde{X}_{\mid _\varepsilon }^\beta }{\partial
x^i}-\widetilde{X}^\alpha \sigma _\alpha ^i \frac{\partial
\mathcal{S}^\beta }{\partial x^i}+\mathcal{S}^\alpha
\widetilde{X}_{\mid _\alpha }^\beta -y^\varepsilon
\widetilde{X}_{\mid _\varepsilon }^\alpha \frac{\partial
\mathcal{S}^\beta }{\partial y^\alpha } =0.
\end{equation*}
and it results, by direct computation, that the components
$\widetilde{X} ^\alpha (x)$ satisfy the equations (23), (24).

\begin{definition}
A section $X\in \Gamma (\mathcal{T}E\backslash \{0\}$ is called
Newtonoid if $J[\mathcal{S},X]_{\mathcal{T}E}=0.$
\end{definition}
In local coordinates we obtain
\begin{equation*}
J[\mathcal{S},X]_{\mathcal{T}E}=\left( \mathcal{S}(X^\alpha
)-Y^\alpha +y^\varepsilon L_{\varepsilon \beta }^\alpha X^\beta
\right) \mathcal{V} _\alpha ,
\end{equation*}
which yields
\begin{equation}
Y^\alpha =\mathcal{S}(X^\alpha )+y^\varepsilon L_{\varepsilon
\beta }^\alpha X^\beta ,\quad X=X^\alpha \mathcal{X}_\alpha
+\left( \mathcal{S}(X^\alpha )+y^\varepsilon L_{\varepsilon \beta
}^\alpha X^\beta \right) \mathcal{V} _\alpha .
\end{equation}
We remark that a section $X\in \Gamma (\mathcal{T}E\backslash
\{0\}$ is a dynamical symmetry if and only if it is a Newtonoid
and satisfies the equation (24). The set of Newtonoid sections
denoted $\frak{X}_{\mathcal{S}}$ is given by
\begin{equation*}
\frak{X}_{\mathcal{S}}=Ker(J\circ \mathcal{L}_{\mathcal{S}})=Im(Id+J\circ
\mathcal{L}_{\mathcal{S}}).
\end{equation*}
In the following we will use the dynamical covariant derivative in
order to find the invariant equations of Newtonoid sections and
dynamical symmetries on Lie algebroids. Let $\mathcal{S}$ be a
semispray, $\mathcal{N}$ an arbitrary nonlinear connection and
$\nabla $ the induced dynamical covariant derivative. We set:

\begin{proposition}
A section $X\in \Gamma (\mathcal{T}E\backslash \{0\})$ is a
Newtonoid if and only if
\begin{equation}
\mathrm{v}(X)=J(\nabla X),
\end{equation}
which locally yields
\begin{equation*}
X=X^\alpha \delta _\alpha +\nabla X^\alpha \mathcal{V}_\alpha,
\end{equation*}
with $\nabla X^\alpha $ given by formula (18).
\end{proposition}

\textbf{Proof}. We know that $J\circ \nabla =J\circ
\mathcal{L}_{\mathcal{S} }+\mathrm{v}$ and it results
$J[\mathcal{S},X]_{\mathcal{T}E}=0$ if and only if
$\mathrm{v}(X)=J(\nabla X).$ In local coordinates we obtain
\begin{eqnarray*}
X &=&X^\alpha \left( \delta _\alpha +\mathcal{N}_\alpha ^\beta
\mathcal{V} _\beta \right) +\left( \mathcal{S}(X^\alpha
)+y^\varepsilon L_{\varepsilon
\beta }^\alpha X^\beta \right) \mathcal{V}_\alpha \\
\ &=&X^\alpha \delta _\alpha +\left( \mathcal{S}(X^\alpha
)+X^\beta ( \mathcal{N}_\beta ^\alpha +y^\varepsilon
L_{\varepsilon \beta }^\alpha
\right) \mathcal{V}_\alpha \\
\ &=&X^\alpha \delta _\alpha +\nabla X^\alpha \mathcal{V}_\alpha.
\end{eqnarray*}
\hfill
\hbox{\rlap{$\sqcap$}$\sqcup$}

\begin{proposition}
A section $X\in \Gamma (\mathcal{T}E\backslash \{0\})$ is a
dynamical symmetry if and only if $X$ is a Newtonoid and
\begin{equation}
\nabla (J\nabla X)+\Phi (X)=0.
\end{equation}
\end{proposition}

\textbf{Proof}. If $X$ is a dynamical symmetry then
$\mathrm{h}[\mathcal{S}
,X]_{\mathcal{T}E}=\mathrm{v}[\mathcal{S},X]_{\mathcal{T}E}=0$ and
composing by $J$ we get $J[\mathcal{S},X]_{\mathcal{T}E}=0$ that
means $X$ is a Newtonoid. Therefore,
$\mathrm{v}[\mathcal{S},X]_{\mathcal{T}E}=\mathrm{v}[
\mathcal{S},\mathrm{v}X]_{\mathcal{T}E}+\mathrm{v}[\mathcal{S},\mathrm{h}X]_{
\mathcal{T}E}=\nabla (\mathrm{v}X)+\Phi (X)$ and using (26) we get
$\nabla (J\nabla X)+\Phi (X)=0.$\hfill
\hbox{\rlap{$\sqcap$}$\sqcup$}\\

For $f\in C^\infty (E)$ and $X\in \Gamma (\mathcal{T}E\backslash
\{0\})$ we define the product
\begin{equation*}
f*X=(Id+J\circ
\mathcal{L}_{\mathcal{S}})(fX)=fX+fJ[\mathcal{S},X]_{\mathcal{
T}E}+\mathcal{S}(f)JX,
\end{equation*}
and remark that a section $X\in \Gamma (\mathcal{T}E\backslash
\{0\})$ is a Newtonoid if and only if
\begin{equation*}
X=X^\alpha (x,y)*\mathcal{X}_\alpha .
\end{equation*}
If $X\in \frak{X}_{\mathcal{S}}$ then
\begin{equation*}
f*X=fX+\mathcal{S}(f)JX.
\end{equation*}
The next result proves that the canonical nonlinear connection can
by determined by symmetry.
\begin{proposition}
Let us consider a semispray $\mathcal{S}$, an arbitrary nonlinear
connection $\mathcal{N}$ and $\nabla $ the dynamical covariant
derivative. The following conditions are equivalent:\\
$i)$ $\nabla$ restricts to $\nabla
:\frak{X}_{\mathcal{S}}\rightarrow \frak{X }_{\mathcal{S}}$
satisfies the Leibnitz rule with respect to the $*$ pro-duct.\\
$ii)$ $\nabla J=0$,\\
$iii)$ $\mathcal{L}_{\mathcal{S}}J+\mathcal{N}=0,$\\
$iv)$ $\mathcal{N}_\alpha ^\beta =\frac 12\left( -\frac{\partial
\mathcal{S} ^\beta }{\partial y^\alpha }+y^\varepsilon L_{\alpha
\varepsilon }^\beta \right)$.
\end{proposition}

\textbf{Proof}. For $ii)\Rightarrow i)$ we consider $X\in
\frak{X}_{\mathcal{ S}}$ and using (26) we have
$\mathrm{v}X=J(\nabla X)$ which leads to $\nabla
(\mathrm{v}X)=\nabla (J\nabla X)$. It results $(\nabla
\mathrm{v})X+\mathrm{v }(\nabla X)=(\nabla J)(\nabla X)+J\nabla
(\nabla X)$ and using the relations $\nabla \mathrm{v}=0$ and
$\nabla J=0$ we obtain $\mathrm{v}(\nabla X)=J\nabla (\nabla X)$
which implies $\nabla X\in \frak{X}_{\mathcal{S}}$ . For $X\in
\frak{X}_{\mathcal{S}}$ we have
\begin{equation*}
\nabla \left( f*X\right) =\nabla (fX+\mathcal{S}(f)JX)=\mathcal{S}
(f)X+f\nabla X+\mathcal{S}^2(f)JX+\mathcal{S}(f)\nabla (JX),
\end{equation*}
\begin{equation*}
\nabla f*X+f*\nabla X=\mathcal{S}(f)X+\mathcal{S}^2(f)JX+f\nabla
X+\mathcal{S }(f)J(\nabla X).
\end{equation*}
But $\nabla (JX)=(\nabla J)X+J(\nabla X)$ and from $\nabla J=0$ we
obtain $ \nabla (JX)=J(\nabla X)$ which leads to $\nabla \left(
f*X\right) =\nabla f*X+f*\nabla X$.

For $i)\Rightarrow ii)$ we consider the set
$\frak{X}_{\mathcal{S}}\cup \Gamma
^{\mathrm{v}}(\mathcal{T}E\backslash \{0\})$ which is a set of
generators for $\Gamma (\mathcal{T}E\backslash \{0\})$. We have
$\nabla J(X)=0$ for $X\in \Gamma
^{\mathrm{v}}(\mathcal{T}E\backslash \{0\})$ and for $X\in
\frak{X}_{\mathcal{S}}$ using $\nabla \left( f*X\right) =\nabla
f*X+f*\nabla X$ it results $\mathcal{S}(f)\nabla
(JX)=\mathcal{S}(f)J(\nabla X),$ which implies
$\mathcal{S}(f)(\nabla J)X=0,$ for an arbitrary function $ f\in
C^\infty (E\backslash \{0\})$. Therefore, $\nabla J=0$ on $
\frak{X}_{\mathcal{S}}$ which ends the proof. The equivalence of
the conditions $ii)$, $iii)$, $iv)$ have been proved in the
Theorem 1. \hfill \hbox{\rlap{$\sqcap$}$\sqcup$}

Next, we consider the dynamical covariant derivative $\nabla$
induced by the semispray $\mathcal{S}$, the canonical nonlinear
connection $\mathcal{N}=- \mathcal{L}_{\mathcal{S}}J$ and find the
invariant equations of dynamical and Lie symmetries.

\begin{proposition}
A section $X\in \Gamma (\mathcal{T}E\backslash \{0\}$ is a
dynamical symmetry if and only if $X$ is a Newtonoid and
\begin{equation}
\nabla ^2JX+\Phi (X)=0,
\end{equation}
which locally yields
\begin{equation*}
\nabla ^2X^\alpha +\mathcal{R}_\beta ^\alpha X^\beta =0.
\end{equation*}
\end{proposition}

\textbf{Proof}. From (20) it results $\nabla J=0$ and using (27)
and (17) we get (28). Next, using (25) and (14) the local
components of the vertical section $\nabla ^2JX+\Phi (X)$ is
$\nabla ^2X^\alpha +\mathcal{R}_\beta ^\alpha X^\beta$. \hfill
\hbox{\rlap{$\sqcap$}$\sqcup$}

\begin{proposition}
A section $\widetilde{X}\in \Gamma (E\backslash \{0\})$ is a Lie
symmetry of $\mathcal{S}$ if and only if
\begin{equation}
\nabla ^2\widetilde{X}^v+\Phi (\widetilde{X}^c)=0
\end{equation}
\end{proposition}

\textbf{Proof}. Using (28) and the relation
$J(\widetilde{X}^c)=\widetilde{X} ^v$ we obtain (29). \hfill
\hbox{\rlap{$\sqcap$}$\sqcup$}\\

Let us consider in the following a regular Lagrangian $L$ on $E$,
the Cartan 1-section $\theta_L$, the symplectic structure
$\omega_L=d^E\theta_L$, the energy function $E_L$ and the induced
canonical semispray $\mathcal{S}$ with the components given by the
relation (9).

\begin{proposition}
If $\widetilde{X}$ is a section on $E$ such that
$\mathcal{L}_{\widetilde{X} ^c}\theta _{L\text{ }}$is closed and
$d^E(\widetilde{X}^cE_L)=0$, then $ \widetilde{X}$ is a Lie
symmetry of the canonical semispray $\mathcal{S}$ induced by $L$.
\end{proposition}

\textbf{Proof}. We have
\begin{eqnarray*}
i_{[\widetilde{X}^c,\mathcal{S}]}\omega _L
&=&\mathcal{L}_{\widetilde{X} ^c}(i_{\mathcal{S}}\omega
_L)-i_{\mathcal{S}}(\mathcal{L}_{\widetilde{X} ^c}\omega
_L)=-\mathcal{L}_{\widetilde{X}^c}d^EE_L-i_{\mathcal{S}}(\mathcal{L
}_{\widetilde{X}^c}d^E\theta _L) \\
&=&-d^E\mathcal{L}_{\widetilde{X}^c}E_L-i_{\mathcal{S}}d^E(\mathcal{L}_{
\widetilde{X}^c}\theta
_L)=-d^E(\widetilde{X}^cE_L)-i_{\mathcal{S}}d^E(
\mathcal{L}_{\widetilde{X}^c}\theta _L)=0.\end{eqnarray*} But
$\omega _L$ is a symplectic structure ($L$ is regular) and we get
$[ \widetilde{X}^c,\mathcal{S}]=0$ which ends the proof. \hfill
\hbox{\rlap{$\sqcap$}$\sqcup$}

\begin{definition}
a) A section $X\in \Gamma (\mathcal{T}E\backslash \{0\})$ is
called a Cartan symmetry of the Lagrangian $L$, if
$\mathcal{L}_X\omega _L=0$ and $\mathcal{L }_XE_L=0$.\\
b) A function $f\in C^\infty (E)$ is a constant of motion (or a
conservation law) for the Lagrangian $L$ if $\mathcal{S}(f)=0$
\end{definition}

\begin{proposition}
The canonical semispray induced by the regular Lagrangian $L$ is a
Cartan symmetry.
\end{proposition}

\textbf{Proof}. Using the relation
$i_{\mathcal{S}}\omega_L=-d^EE_L$ and the skew symmetry of the
symplectic 2-section $\omega _L$ we obtain
\[
0=i_{\mathcal{S}}\omega
_L(\mathcal{S})=-d^EE_L(\mathcal{S})=-\mathcal{S}(E_L)=-
\mathcal{L}_SE_L.
\]
Also, from $d^E\omega _L=0$ we get
\[
\mathcal{L}_{\mathcal{S}}\omega _L=d^Ei_{\mathcal{S}}\omega
_L+i_{\mathcal{S} }d^E\omega _L=-d^E(d^EE_L)=0,
\]
and it results that the semispray $\mathcal{S}$ is a Cartan
symmetry.\hfill \hbox{\rlap{$\sqcap$}$\sqcup$}

\begin{proposition}
A Cartan symmetry $X$ of the Lagrangian $L$ is a dynamical
symmetry for the canonical semispray $\mathcal{S}$.
\end{proposition}

\textbf{Proof}. From the symplectic equation $i_S\omega
_L=-d^EE_L$, applying the Lie derivative in both sides, we obtain
\begin{equation*}
\mathcal{L}_X(i_S\omega
_L)=-\mathcal{L}_Xd^EE_L=-d^E\mathcal{L}_XE_L=0.
\end{equation*}
Also, using the formula $i_{[X,Y]_{\mathcal{T}E}}=\mathcal{L}_X\circ
i_Y-i_Y\circ \mathcal{L}_X$ it results
\begin{equation*}
\mathcal{L}_X(i_{\mathcal{S}}\omega
_L)=-i_{[\mathcal{S},X]_{\mathcal{T} E}}\omega
_L+i_{\mathcal{S}}\mathcal{L}_X\omega _L=-i_{[\mathcal{S},X]_{
\mathcal{T}E}}\omega _L
\end{equation*}
which yields
\begin{equation}
i_{[\mathcal{S},X]_{\mathcal{T}E}}\omega _L=0.
\end{equation}
But $\omega _L$ is a symplectic 2-section and we conclude that
$[\mathcal{S},X]_{ \mathcal{T}E}=0$, so $X$ is a
dynamical symmetry.\hfill \hbox{\rlap{$\sqcap$}$\sqcup$}\\
Since Lie and exterior derivatives commute, we obtain
\[
d^E\mathcal{L}_X\theta _L=\mathcal{L}_Xd^E\theta
_L=\mathcal{L}_X\omega _L=0,
\]
It results that, for a Cartan symmetry, the 1-section
$\mathcal{L}_X\theta _L$ is a closed 1-section.
\begin{definition}
A Cartan symmetry $X$ is said to be an exact Cartan symmetry if
the 1-section $\mathcal{L}_X\theta _L$ is exact.
\end{definition}
The next result proves that there is a one to one correspondence
between exact Cartan symmetries and conservation laws. Also, if
$X$ is an exact Cartan symmetry, then there is a function $f\in
C^\infty (E)$ such that $\mathcal{L}_X\theta _L=d^{E}f$.
\begin{proposition}
If $X$ is an exact Cartan symmetry, then $f-\theta_L(X)$ is a
conservation law for the Lagrangian $L$. Conversely, if $f\in
C^\infty (E)$ is a conservation law for $L$, then $X\in \Gamma
(\mathcal{T}E\backslash \{0\})$ the unique solution of the
equation $i_X\omega_L=-d^Ef$ is an exact Cartan symmetry.
\end{proposition}
\textbf{Proof}. We have $\mathcal{S}(f-\theta _L(X)) =d^E(f-\theta
_L(X))(\mathcal{S})=\left( \mathcal{L}_X\theta _L-d^Ei_X(\theta
_L)\right) (\mathcal{S})=i_Xd^E\theta _L(\mathcal{S}) =i_X\omega
_L(\mathcal{S})=-i_{\mathcal{S}}\omega _L(X)=d^EE_L(X)=0,$ and it
results that $f-\theta _L(X)$ is a conservation law for the
dynamics associated to the regular Lagrangian $L$. Conversely, if
$X$ is the solution of the equation $i_X\omega_L=-d^Ef$ then
$\mathcal{L}_X{\theta_L}=i_X{\omega_L}$ is an exact 1-section.
Consequently, $0=d^E\mathcal{L}_X\theta _L=\mathcal{L}_Xd^E\theta
_L=\mathcal{L}_X\omega _L.$ Also, $f$ is a conservation law, and
we have $0=\mathcal{S}(f)=d^Ef(\mathcal{S})=-i_X\omega
_L(\mathcal{S})=i_S\omega _L(X)=-d^EE_L(X)=-X(E_L).$ Therefore, we
obtain $\mathcal{L}_XE_L=0$ and $X$ is an exact Cartan symmetry.
\hfill\hbox{\rlap{$\sqcap$}$\sqcup$} \\We have to mention that the
Noether type theorems for Lagrangian systems on Lie algebroids are
studied in \cite{Ca, Ma2} and Jacobi sections for second order
differential equations on Lie algebroids are investigated in
\cite{Ca1}. \quad

\subsection{Example}

Next, we consider an example from optimal control theory and prove
that the framework of Lie algebroids is more useful that the
tangent bundle in order to calculate some symmetries of the
dynamics induced by a Lagrangian function. Let us consider the
following distributional system in $\Bbb{R}^3$ (driftless control
affine system) \cite{Po2'}:
\[
\left\{
\begin{array}{l}
\dot x^1=u^1+u^2x^1 \\
\dot x^2=u^2x^2 \\
\dot x^3=u^2
\end{array}
\right.
\]
Let $x_0$ and $x_1$ be two points in $\Bbb{R}^3$. An optimal
control problem consists of finding the trajectories of our
control system which connect $x_0 $ and $x_1$ and minimizing the
Lagrangian
\[
{\min }\int_0^T\mathcal{L}(u(t))dt,\ \mathcal{L} (u)=\frac
12\left( (u^1)^2+(u^2)^2\right) ,\quad x(0)=x_0,\ x(T)=x_1,
\]
where $\dot x^i=\frac{dx^i}{dt}$ and $u^1,u^2$ are control
variables. From the system of differential equations we obtain
$u^2=\dot x^3$, $u^1=\dot x^1-\dot x^3x^1$. The Lagrangian
function on the tangent bundle $T\Bbb{R}^3$ has the form
\[
\mathcal{L}=\frac 12\left( \dot x^1-\dot x^3x^1)^2+(\dot
x^3)^2\right) ,
\]
with the constraint
\[
\dot x^2=\dot x^3x^2.
\]
Then, using the Lagrange multiplier $\lambda =\lambda (t)$, we
obtain the total Lagrangian (including the constraints) given by
\[
L(x,\dot x)=\mathcal{L}(x,\dot x)+\lambda \left( \dot x^2-\dot
x^3x^2\right) =\frac 12\left( (\dot x^1-\dot x^3x^1)^2+(\dot
x^3)^2\right) +\lambda \left( \dot x^2-\dot x^3x^2\right) .
\]
We observe that the Hessian matrix of $L$ is singular, and $L$ is
a degenerate Lagrangian (not regular). The corresponding
Euler-Lagrange equations lead to a  complicated system of second
order differential equations. Moreover, because the Lagrangian is
not regular, we cannot obtain the explicit coefficients of the
semispray $\mathcal{S}$ from the equation $i_{\mathcal{S}}\omega
_L=-dE_L$ and it is difficult to study the symmetries of SODE in
this case.\\
For this reason, we will use a different approach, considering the
framework of Lie algebroids. The system can be written in the next
form
\begin{eqnarray*}
\left. \dot x=u^1X_1+u^2X_2,\quad x=\left(
\begin{array}{c}
x^1 \\
x^2 \\
x^3
\end{array}
\right) \in \Bbb{R}^3,\ X_1=\left(
\begin{array}{c}
1 \\
0 \\
0
\end{array}
\right) ,\ X_2=\left(
\begin{array}{c}
x^1 \\
x^2 \\
1
\end{array}
\right) .\right.
\end{eqnarray*}
The associated distribution $\Delta =span\{X_1,X_2\}$ has the
constant rank $2$ and is holonomic, because
\[
X_1=\frac \partial {\partial x^1},\quad X_2=x^1\frac \partial
{\partial x^1}+x^2\frac \partial {\partial x^2}+\frac \partial
{\partial x^3},\quad [X_1,X_2]=X_1.
\]
From the Frobenius theorem, the distribution $\Delta $ is
integrable, it determines a foliation on $T\Bbb{R}^3$ and two
points can be joined by a optimal trajectory if and only if they
are situated on the same leaf (see \cite{Po2'}). In order to apply
the theory of Lie algebroids, we consider the Lie algebroid being
just the distribution, $E=\Delta $ and the anchor $\sigma
:E\rightarrow T\Bbb{R}^3$ is the inclusion, with the components
\[
\sigma _\alpha ^i=\left(
\begin{array}{cc}
1 & x^1 \\
0 & x^2 \\
0 & 1
\end{array}
\right) .
\]
From the relation
\[
\lbrack X_\alpha ,X_\beta ]=L_{\alpha \beta }^\gamma X_\gamma
,\quad \alpha ,\beta ,\gamma =1,2,
\]
we obtain the non-zero structure functions
\[
L_{12}^1=1,\ L_{21}^1=-1.
\]
The components of the semispray from (9) are given by
\[
\mathcal{S}^1=-u^1u^2,\ \mathcal{S}^2=(u^1)^2.
\]
The functions $\mathcal{S}^\alpha $ are homogeneous of degree 2 in
$u$ and it results that $\mathcal{S}$ is a spray. By
straightforward computation we obtain the expression of the
canonical spray induced by $\mathcal{L}$
\[
\mathcal{S}(x,u)=(u^1+u^2x^1)\frac \partial {\partial
x^1}+u^2x^2\frac
\partial {\partial x^2}+u^2\frac \partial {\partial x^3}-u^1u^2\frac
\partial {\partial u^1}+(u^1)^2\frac \partial {\partial u^2}.
\]
From the Proposition 17 it results that $\mathcal{S}(x,u)$ is a
Cartan symmetry of the dynamics associated to the regular
Lagrangian $\mathcal{L}$ on Lie algebroids.\\
The  coefficients of the canonical nonlinear connection
$\mathcal{N}=- \mathcal{L}_{\mathcal{S}}J$ are given by
\[
\mathcal{N}_1^1=u^2,\ \mathcal{N}_2^1=0,\ \mathcal{N}_1^2=u^1,\
\mathcal{N}_2^2=0,
\]
and the components of the Jacobi endomorphism from (21) have the
form
\[
\mathcal{R}_1^1=-(u^2)^2,\ \mathcal{R}_1^2=-u^1u^2,\ \mathcal{R}
_2^1=u^1u^2,\ \mathcal{R}_2^2=(u^1)^2.
\]
Also, the non-zero coefficients of the curvature from (11) of
$\mathcal{N}$ are
\[
\mathcal{R}_{12}^1=u^2,\ \mathcal{R}_{12}^2=u^1,\
\mathcal{R}_{21}^1=-u^2,\ \mathcal{R}_{21}^2=-u^1,
\]
and we obtain that the Jacobi endomorphism is the contraction with
$\mathcal{S}$ of the curvature of $\mathcal{N}$, or locally
$\mathcal{R}_\beta ^\alpha =\mathcal{R}_{\varepsilon \beta
}^\alpha u^\varepsilon $.\\
The Euler-Lagrange equations on Lie algebroids given by (see
\cite{We})
\[
\frac{dx^i}{dt}=\sigma _\alpha ^iu^\alpha ,\quad \frac d{dt}\left(
\frac{
\partial \mathcal{L}}{\partial u^\alpha }\right) =\sigma _\alpha ^i\frac{
\partial \mathcal{L}}{\partial x^i}-L_{\alpha \beta }^\varepsilon u^\beta
\frac{\partial \mathcal{L}}{\partial u^\varepsilon },
\]
lead to the following differential equations
\[
\dot u^1=-u^1u^2,\quad \dot u^2=(u^1)^2,
\]
which can be written in the form
\[
\frac{dx^i}{dt}=\sigma _\alpha ^iu^\alpha ,\quad \frac{du^\alpha
}{dt}= \mathcal{S}^\alpha (x,u),
\]
and give the integral curves of $\mathcal{S}$. The Cartan
1-section $\theta _{\mathcal{L}}$ has the form
\[
\theta _{\mathcal{L}}=u^1dx^1+u^2(x^1dx^1+x^2dx^2+dx^3),
\]
and the symplectic structure is $\omega _{\mathcal{L}}=d^E\theta
_{\mathcal{L}}$. The energy of Lagragian $\mathcal{L}$ is
\[
E_{\mathcal{L}}=\frac 12\left( (u^1)^2+(u^2)^2\right).
\]
For the optimal
solution of the control system
(using the framework of Lie algebroids) see \cite{Po2'}.\\
 \quad \\ \textbf{Conclusions}. The main purpose of this
paper is to study the symmetries of SODE on Lie algebroids and
relations between them, using the dynamical covariant derivative
and Jacobi endomorphism. The existence of a semispray
$\mathcal{S}$ together with an arbitrary nonlinear connection $
\mathcal{N}$ define a dynamical covariant derivative and the
Jacobi endomorphism. Let us remark that at this point we do not
have any relation between $\mathcal{S}$ and the nonlinear
connection $\mathcal{N}.$ This will be given considering the
compatibility condition between the dynamical covariant derivative
and the tangent structure, $\nabla J=0$, which fix the canonical
nonlinear connection $\mathcal{N}=-\mathcal{L}_{ \mathcal{S}}J$.
This canonical nonlinear connection depends only on semispray. In
this case we have the decomposition $\nabla =\mathcal{L}_{
\mathcal{S}}+\Bbb{F}+\mathcal{J}-\Phi $ which can be compared with
the tangent case from \cite{Bu3, Ma}. Also, in the case of
homogeneous SODE (spray), the dynamical covariant derivative
coincides with Berwald linear connection and the Jacobi
endomorphism is the contraction with $\mathcal{S}$ of the
curvature of the nonlinear connection. We study the dynamical
symmetry, Lie symmetry, Newtonoid section and Cartan symmetry on
Lie algebroids and find their invariant equations with the help of
dynamical covariant derivative and Jacobi endomorphism. Finally,
we give an example from optimal control theory which proves that
the framework of Lie algebroids is more useful than the tangent
bundle in order to find the symmetries of the dynamics induced by
a Lagrangian function. For further developments one can study the
symmetries using the
$k$-symplectic formalism on Lie algebroids given in \cite{Le2}.\\
\textbf{Acknowledgments}. The author wishes to express his thanks
to the referees for useful comments and suggestions concerning
this paper.

Author's address: \\University of Craiova,\\Dept. of Statistics and Informatics,\\
13, Al. I. Cuza, st., Craiova 200585, Romania\\e-mail:
liviupopescu@central.ucv.ro; liviunew@yahoo.com


\begin{thebibliography}{99}
\bibitem{Ab}  R. Abraham, J. Marsden, \textit{Foundation of Mechanics},
Benjamin, New-York, 1978.

\bibitem{Ar2}  C.M. Arcu\c s, \textit{Mechanical systems in the generalized
Lie algebroids framework,} Int. J. Geom. Methods Mod. Phys., 11, 1450023
(2014).

\bibitem{Bar} M. Barbero-Li\~n\'an, M. Farr\'e Puiggal\'\i, D. Mart\'\i n de Diego,
\textit{Inverse problem for Lagrangian systems on Lie algebroids
and applications to reduction by symmetries}, Monatsh Math., vol.
180, 4 (2016) 665-691.

\bibitem{Bua}  L. Bua, I. Bucataru, M. Salgado, \textit{Symmetries,
Newtonoid vector fields and conservation laws on the Lagrangian k-symplectic
formalism}, Reviews in Mathematical Physics, 24, 1250030 (2012).

\bibitem{Bu2}  I. Buc\u ataru, M.F. Dahl, \textit{Semi-bazic 1-form and
Helmholtz conditions for the inverse problem of the calculus of variations},
J. Geom. Mechanics, 1, no.2 (2009), 159-180.

\bibitem{Bu3}  I. Buc\u ataru, O. Constantinescu, M.F. Dahl, \textit{A
geometric setting for systems of ordinary differential equations,} Int. J.
Geom. Methods Mod. Phys., 08, 12019 (2011).

\bibitem{Ca0} J. F. Cari\~nena, E. Mart\'\i nez, \textit{Generalized Jacobi equation
and inverse problem in classical mechanics}, in Group Theoretical
Methods in Physics, V. V. Dodonov and V. I. Manko, eds.,
Proceedings of the 18th International Colloquim 1990, Moskow, Vol.
2 (Nova Science Publishers, New York, 1991).

\bibitem{Ca}  J.F. Cari\~nena, M. Rodr\'\i guez-Olmos, \textit{Gauge
equivalence and conserved quantities for Lagrangian systems on Lie
algebroids} , J. Phys. A: Math. Theor. 42 (2009), 265209.

\bibitem{Ca1} J.F. Cari\~nena, I. Gheorghiu, E. Mart\'\i nez, \textit{Jacobi
fields for second-order differential equations on Lie algebroids},
Dynamical Systems, Differential Equations and Applications AIMS
Proceedings, 2015, 213-222.

\bibitem{Co}  J. Cortes, E. Mart\'\i nez,\textit{\ Mechanical control systems on
Lie algebroids}, IMA J. Math. Control Inform. 21 (2004) 457-492.

\bibitem{Cr1}  M.Crampin, \textit{Tangent bundle geometry for Lagrangian
dynamics}, J. Phys. A: Math. Gen. 16 (1983), 3755-3772.

\bibitem{Cr2}  M. Crampin, E. Mart\'\i nez, W. Sarlet, \textit{Linear
connections for system of second-order ordinary differential equations},
Ann. Inst. Henry Poincare, 65 no.2 (1996), 223-249.

\bibitem{Fe}  R. L. Fernandes, \textit{Lie algebroids, holonomy and
characteristic classes}, Advances in Mathematics, 170 (2002) 119-179.

\bibitem{Fr}  A. Fr\"olicher, A. Nijenhuis, \textit{Theory of vector-valued
differential forms}, Nederl. Akad. Wetensch. Proc. Ser. A. 59 (1956),
338-359.

\bibitem{Ga}  X. Gr\`acia, J.M. Pons, \textit{Symmetries and infinitesimal
symmetries of singular differential equations}, J. Phys. A: Math. Gen. 35
(2002), 5059-5077.

\bibitem{Gu1}  J. Grabowski, P. Urbanski, \textit{Tangent and cotangent lift
and graded Lie algebra associated with Lie algebroids}, Ann. Global Anal.
Geom. 15, (1997), 447-486.

\bibitem{Gu2}  J. Grabowski, P. Urbanski, \textit{Lie algebroids and
Poisson-Nijenhuis structures}, Rep. Math. Phys., 40 (1997), 195-208.

\bibitem{Gr1}  J. Grifone, \textit{Structure presque tangente et connections
I}, Ann. Inst. Fourier, 22, no.1 (1972), 287-334.

\bibitem{Gr2}  J. Grifone, Z. Muzsnay, \textit{Variational principle for
second-order differential equations}, World Scientific, 2000.

\bibitem{Hi}  P. J. Higgins, K. Mackenzie, \textit{Algebraic constructions
in the category of Lie algebroids}, Journal of Algebra, 129 (1990), 194-230.

\bibitem{Je} M. Jerie, G. Prince, \textit{Jacobi fields and linear connections
for arbitrary second-order ODEs}, J. Geom. Physics, 43 (2002)
351-370.

\bibitem{Le0}  M. de Le\'on, D. Mart\'\i n de Diego, \textit{Symmetries and
constants of the motion for singular Lagrangian systems,} Int. J. Theor.
Phys. 35 (5) (1996), 975-1011.

\bibitem{Le}  M. de Le\'on, J. C. Marrero, E. Mart\'\i nez, \textit{Lagrangian
submanifolds and dynamics on Lie algebroids}, J. Phys. A: Math. Gen. 38
(2005), 241--308.

\bibitem{Le1}  M. De Le\'on, D. Mart\'\i n de Diego, A. Santamaria-Merino,
\textit{Symmetries in classical field theories, } Int. J. Geom. Meth. Mod.
Phys., 5 (2004), 651-710\textit{.}

\bibitem{Le2}  M. de Le\'on, D. Mart\'\i n de Diego, M. Salgado, S.
Vilari\~no, \textit{k-symplectic formalism on Lie algebroids.} J. Phys. A:
Math. Theor. 42 (2009) 385209

\bibitem{Li}  P. Libermann, \textit{Lie algebroids and mechanics}, Arch.
Math. Brno, 32 (1996), 147-162.

\bibitem{Mk1}  K. Mackenzie, \textit{Lie groupoids and Lie algebroids in
differential geometry}, London Mathematical Society Lecture Note Series 124,
1987.

\bibitem{Mk2}  K. Mackenzie, \textit{General theory of Lie groupoids and Lie
algebroids}, London Mathematical Society, Cambridge 213, 2005.

\bibitem{Mar} G. Marmo, N. Mukunda, \textit{Symmetries and constant of the
motion in the Lagrangian formalism on $TQ$: beyond point
transformations}, Nuovo Cim. B, 92 (1986) 1-12.

\bibitem{Ma}  E. Mart\'\i nez, J. F. Cari\~nema, W. Sarlet, \textit{
Derivations of differential forms along the tangent bundle projection II},
Diff. Geom. Appl. 3, no.1 (1993), 1-29.

\bibitem{Ma2}  E. Mart\'\i nez, \textit{Lagrangian mechanics on Lie algebroids},
Acta Appl. Math. 67 (2001), 295--320.

\bibitem{Ma4}  E. Mart\'\i nez, T. Mestdag, W. Sarlet, \textit{Lie algebroid
structures and Lagrangian systems on affine bundles}, J. Geom. Phys. 44, no.
1, (2002), 70--95.

\bibitem{Me}  T. Mestdag, B. Langerock, \textit{A Lie algebroid framework
for non-holonomic systems}, J. Phys. A: Math. Gen., 38, 5, 1097 (2005).

\bibitem{Pe}  E. Peyghan, \textit{Berwald-type and Yano-type connections on Lie algebroids},
Int. J. Geom. Methods Mod. Phys., vol. 12, no. 10, 1550125 (2015).

\bibitem{Po1}  L. Popescu, \textit{The geometry of Lie algebroids and
applications to optimal control}, Annals. Univ. Al. I. Cuza, Iasi, series I,
Math., LI (2005), 155-170.

\bibitem{Po2}  L. Popescu, \textit{Geometrical structures on Lie algebroids},
 Publ. Math. Debrecen, 72, 1-2 (2008), 95-109.

\bibitem{Po2'}  L. Popescu, \textit{Lie algebroids framework for
distributional systems,} Annals Univ. Al. I. Cuza, Iasi, series I, Math.,
55, 2 (2009), 373-390.

\bibitem{Po3}  L. Popescu, \textit{Metric nonlinear connections on Lie
algebroids}, Balkan J. Geom. Appl. 16, no. 1, (2011), 111-121.

\bibitem{Po4}  L. Popescu, \textit{Dual structures on the prolongations of a
Lie algebroid}, Annals Univ. Al. I. Cuza, Iasi, series I, Math., 59, 2
(2013), 357-374.

\bibitem{Pr1} G. Prince, \textit{Toward a classification of
dynamical symmetries in classical mechanics}, Bull. Austral. Math.
Soc., vol. 27 (1983), 53-71.

\bibitem{Pr2} G. Prince, \textit{A complete classification of
dynamical symmetries in classical mechanics}, Bull. Austral. Math.
Soc., vol. 32 (1985), 299-308.

\bibitem{Ro}  N. Rom\'an-Roy, M. Salgado, S. Vilari\~no, \textit{Symmetries
and conservation laws in G\"unter k-symplectic formalism of field theory},
Reviews in Mathematical Physics, 19 (10) (2007), 1117-1147.

\bibitem{Sa}  W. Sarlet, \textit{Linear connections along the tangent bundle
projection}, In: Variations, Geometry and Physics, Nova Science Publishers
(2008).

\bibitem{Sz2}  J. Szilasi, \textit{A setting for spray and Finsler geometry,}
In: Handbook of Finsler Geometry, Kluwer Acad. Publ., 2 (2003), 1183-1426.

\bibitem{We}  A. Weinstein, \textit{Lagrangian mechanics and groupoids},
Fields Inst. Comm. 7 (1996), 206-231.
\end{thebibliography}
\end{document}